\documentclass[a4paper,11pt]{article}

\usepackage{amsmath,amssymb,subcaption,graphicx}
\usepackage[margin=2.5cm]{geometry}

\usepackage{multicol}

\begin{document}

\title{Mathematical modelling of the vitamin C clock reaction: a study of two kinetic regimes}

\author{
A. A. Alsaleh$^{1}$, D. J. Smith\thanks{Correspondence: \texttt{d.j.smith@bham.ac.uk}}$^{\ 1,2}$ and S. Jabbari$^{1}$\\
$^{1}$School of Mathematics, 
$^{2}$Centre for Systems Modelling and\\Quantitative Biomedicine, University of Birmingham, Edgbaston,\\ Birmingham, B15 2TT, United Kingdom}

\maketitle

\begin{abstract}
Chemically reacting systems exhibiting a repeatable delay period before a visible and sudden change are referred to as \emph{clock reactions}; they have a long history in education and provide an idealisation of various biochemical and industrial processes. We focus on a purely substrate-depletive clock reaction utilising vitamin C, hydrogen peroxide, iodine and starch. Building on a recent study of a simplified two-reaction model under high hydrogen peroxide concentrations, we develop a more detailed model which breaks the slow reaction into two steps, one of which is rate-limiting unless hydrogen peroxide levels are very high. Through asymptotic analysis, this model enables the effect of hydrogen peroxide concentration to be elucidated in a principled way, resolving an apparent discrepancy with earlier literature regarding the order of the slow reaction kinetics. The model is analysed in moderate- and high-hydrogen peroxide regimes, providing approximate solutions and expressions for the switchover time which take into account hydrogen peroxide concentration. The solutions are validated through simultaneously fitting the same set of parameters to several experimental series, then testing on independent experiments across widely varying hydrogen peroxide concentration. The study thereby presents and further develops a validated mechanistic understanding of a paradigm chemical kinetics system.
\end{abstract}

\section{Introduction}
Clock reactions are characterised by a defined and predictable induction period followed by a sudden and typically visible change in reactant concentrations. The study of these reactions dates back at least to the work of Landolt on the sulfite/iodate reaction in the 1880s \cite{landolt1886ueber}. These systems have long been used in chemistry education \cite{conway1940modified}, have industrial applications \cite{billingham1993simple,billingham1994kinetics} and alongside experiment have been studied through mathematical and computational models \cite{chien1948kinetic,anderson1976computer,billingham1992mathematical,billingham1993mathematical,kerr2019mathematical}. Clock reactions are of current interest through recent applications as diverse as a chemical clock car activity used in an educational setting \cite{parra2020enhancing}, the evaluation of 3D printed mixing devices \cite{farris2020mix} and the determination of microconcentrations of the potentially toxic dye indigo carmine \cite{pagnacco2022indigo}.
For further details, see the reviews \cite{horvath2015classification,molla2022chemical}.

In this article we will focus on a specific version of the vitamin C clock reaction, described by Wright \cite{wright2002tick,wright2002vitamin} with the aim of providing a safe version of the experiment using only the `household chemicals' of iodine, hydrogen peroxide, starch and ascorbic acid (vitamin C). The system involves two opposing reaction processes: a `fast reaction' which consumes vitamin C (\(C_{6}H_{8}O_{6}\)) and converts molecular iodine (\(I_2\)) to iodide ions (\(I^{-}\)),
\begin{equation}
I_{2}(aq)+C_{6}H_{8}O_{6}(aq)\rightarrow2H^{+}(aq)+2I^{-}(aq)+ C_{6}H_{6}O_{6}(aq).\label{eq:reaction1}
\end{equation}
and a `slow reaction' which consumes hydrogen peroxide (\(H_2O_2\)) and converts iodide ions to molecular iodine,
\begin{equation}
2H^{+}(aq)+2I^{-}(aq)+ H_{2}O_{2}(aq)\rightarrow I_{2}(aq)+2H_{2}O(l).\label{eq:reaction2}
\end{equation}
As will be discussed below, the overall reaction \eqref{eq:reaction2} encapsulates the effect of a rate-determining reaction between iodide ions and hydrogen peroxide to produce hypoiodous acid, and a second faster step involving the combination of hypoiodous acid with a second iodide ion to produce molecular iodine \cite{copper1998kinetics}.

Following initial mixing of iodine, hydrogen peroxide, starch and vitamin C, the fast reaction initially dominates, ensuring that the iodide concentration is much higher than the iodine concentration; in the presence of starch the solution appears white. The solution remains in this state until the vitamin C concentration is nearly depleted, at which point the slow reaction then dominates, converting iodide ions to molecular iodine. In the presence of starch, the colour of the solution darkens so that it appears blue. The time interval at which this colour change occurs will be referred to as the `switchover time' \(t_{sw}\); this interval is repeatable, meaning that the system fits Lente's strict definition of a clock reaction \cite{lente2007}. Within the terminology of Horv\'{a}th and Nagyp\'{a}l \cite{horvath2015classification}, the system is \emph{purely substrate depletive}, by contrast with systems which involve autocatalysis or other clock mechanisms. The subject of this article is predicting on the basis of a mathematical model how \(t_{sw}\) depends on the initial concentrations of the reactants. For other clock reactions involving vitamin C, see references \cite{burgess2014kinetics,farris2020mix}.

The vitamin C clock reaction has been modelled mathematically and tested experimentally \cite{kerr2019mathematical}; this work focused on a regime in which hydrogen peroxide levels are greatly in excess compared with other reactants, which enabled a simplified model involving only two variables, iodine and vitamin C. By contrast, a recent pedagogical application \cite{parra2020enhancing} worked in a regime in which hydrogen peroxide levels were comparable to the other reactants. In addition to requiring the inclusion of an additional variable to take account of the consumption of hydrogen peroxide, this regime makes a qualitatitive change to the apparent kinetics of the slow reaction, from being effectively quadratic in the high hydrogen peroxide case to effectively linear with moderate hydrogen peroxide. The question of whether the reaction should be modelled as linear or quadratic arose in the peer review of reference \cite{kerr2019mathematical}, which proved a significant point of contention due to previous experimental work showing that this reaction is linear \cite{liebhafsky1933kinetics,copper1998kinetics,sattsangi2011microscale}. As discussed in ref.\ \cite{kerr2019mathematical} and the associated open peer review, quadratic reaction kinetics were found to lead to a formula for the switchover time which much better fitted the experimental data of \cite{kerr2019mathematical} than a linear version. It was hypothesised that the reason for this discrepancy was the major difference in hydrogen peroxide concentration used by Kerr et al.\ \cite{kerr2019mathematical} relative to these earlier papers. The hypothesis was that the hypoiodous acid generating step is rate-limiting only when hydrogen peroxide levels are not in great excess.

In this article we will develop a mathematical model which accounts for the details of the production and conversion of hypoiodous acid within the slow reaction, thereby providing a unifying framework for both moderate and high hydrogen peroxide concentrations. By a similar approach to Billingham and Needham \cite{billingham1992mathematical} and Kerr et al.\ \cite{kerr2019mathematical} we will analyse the system through matched asymptotic expansions, exploiting the disparity of slow and fast reaction rates, to develop approximate analytical expressions for the timecourse of the reactant concentrations from initial mixing to final equilibrium. This analysis will lead to approximate formulae for the switchover time in moderate and high hydrogen peroxide regimes in terms of the initial concentrations, one reaction rate parameter, and a parameter relating to the initial iodide/iodine ratio (thus also identifying the dominant reactions and concentrations in the system that drive the switchover time). The mathematical models will be tested through their ability to fit experiment data with the same shared parameter values.

\section{Model formulation}
As in ref.~\cite{kerr2019mathematical}, the fast reaction~\eqref{eq:reaction1} will be expressed in terms of model variables as,
\begin{equation}
I+C\xrightarrow{k_1IC}2D, \label{eq:fastReaction}
\end{equation}
where \(I(t)\) is concentration of molecular iodine, \(C(t)\) is concentration of vitamin C, \(D(t)\) is concentration of iodide and \(k_1\) is the associated reaction rate, as in ref.~\cite{kerr2019mathematical}.

The departure point for the present work is in separating the slow reaction~\eqref{eq:reaction2} into two steps,
\begin{subequations}
\begin{align}
    I^{-}+H_2O_2 &\xrightarrow{k_2} HOI+ OH^-,\label{eq:new1}\\
    I^{-}+HOI & \xrightarrow{k_3} I_2+ OH^-,\label{eq:new2}
\end{align}
\end{subequations}
where the first step is the production of hypoiodous acid \(HOI\), and the second step completes the overall conversion of two iodide ions to one molecule of iodine. The first step is taken to be rate-limiting under conditions of moderate hydrogen peroxide concentration, which can be expressed mathematically as \(k_2/k_3 \ll 1\).

Introducing the variables \(P(t)\) for hydrogen peroxide and \(Q(t)\) for hypoiodous acid, these reactions are simplified as,
\begin{subequations}
\begin{align}
D+P&\xrightarrow{k_2DP}Q, \label{eq:slowReaction1} \\
D+Q&\xrightarrow{k_3DQ}I, \label{eq:slowReaction2}
\end{align}
\end{subequations}
(note that \(OH^-\) is not considered in the model as it has no influence on the other variables).

Finally, there exists a reverse pathway by which hypoiodous acid is reduced by hydrogen peroxide to yield iodide \cite{shin2020kinetics},
\begin{equation}\label{eq:inverse}
    HOI+H_2O_2\rightarrow I^{-},
\end{equation}
or in model variables,
\begin{equation}\label{eq:slowReverse}
Q+P\xrightarrow{k_4QP}D.
\end{equation}

Applying the law of mass action to equations~\eqref{eq:fastReaction}, \eqref{eq:slowReaction1}, \eqref{eq:slowReaction2} and \eqref{eq:slowReverse}
yields the system of ordinary differential equations,
\begin{subequations}\label{eq:M-HP-nondim}
\begin{align}
\frac{d D}{dt}&= -k_2 DP - k_3 DQ + k_4 QP + 2k_1 IC,\\
\frac{d P}{dt}&=  -k_2 DP - k_4 QP,\\
\frac{d Q}{dt}&=  k_2 DP - k_3 DQ - k_4 QP,\\
\frac{d C}{dt}&=  -k_1 IC,\\
\frac{d I}{dt}&=  k_3 DQ-k_1 IC,
\end{align}
\end{subequations}
with initial conditions expressed as
\begin{align*}
D(0)=d_0, \quad P(0)=p_{0},\quad Q(0)=q_{0}, \quad C(0)=c_{0},\quad I(0)=\iota_{0}.
\end{align*}

In summary, the model contains five dependent variables \(D(t), P(t), Q(t), C(t), I(t)\) and has four free parameters \(k_1,k_2,k_3,k_4\) for the reaction rates. The initial values \(p_0\) and \(c_0\) are straightforward to control experimentally, as is the total atomic iodine \(n_0:=2\iota_0+d_0\); it will be assumed below that the initial concentration of hypoiodous acid \(q_0=0\). We will therefore consider \(\phi:=\iota_0/n_0\) to be a fifth free parameter.

Note that conservation of atomic iodine, \(D(t)+Q(t)+2I(t)\) is constant and equal to its initial value \(n_0\). The model can therefore be reduced to
\begin{subequations}\label{eq:update5fullmodel}
\begin{align}
\frac{d P}{dt}&=  -k_2\left(N_{0}-Q-2I\right)P-k_4 QP,\\
\frac{d Q}{dt}&=  k_2\left(N_{0}-Q-2I\right)P-k_3\left(N_{0}-Q-2I\right)Q-k_4 QP,\\
\frac{d C}{dt}&=  -k_1 IC,\\
\frac{d I}{dt}&=  k_3\left(N_{0}-Q-2I\right)Q - k_1 IC.
\end{align}
\end{subequations}

We will make the following assumptions regarding the relative rates of the reactions in the system. The fast reaction is assumed to take place one order of magnitude faster than the fastest part of the slow reaction. The rate limiting step \eqref{eq:slowReaction1} is assumed to be one further order of magnitude slower, as is the reverse reaction \eqref{eq:slowReverse}; schematically,
\begin{equation}
    k_2, k_4 \ll k_3 \ll k_1.   
\end{equation}
More precisely, nondimensionalising with scalings \begin{equation}
P=p_{0} P^{*},\quad Q=n_{0} Q^{*},\quad C=c_{0} C^{*},\quad I=n_{0} I^{*},\quad t=\frac{t^{*}}{k_{1}c_{0}},\label{eq:scalings}
\end{equation}
where $^*$ denotes a dimensionless quantity, yields the dimensionless system,
\begin{subequations}\label{eq:5fullmodelnond}
\begin{align}
\frac{d P^*}{dt^*}&=-\epsilon^2 \sigma \left(1-Q^*-2I^*\right)P^*-\epsilon^2 \beta \sigma Q^*P^*,\\
\frac{d Q^*}{dt^*}&=\epsilon^2 \rho \left(1-Q^*-2I^*\right)P^*-\epsilon \gamma \sigma \left(1-Q^*-2I^*\right)Q^*-\epsilon^2 \beta \rho Q^*P^*,\\
\frac{d C^*}{dt^*}&=-\sigma I^*C^*,\\
\frac{d I^*}{dt^*}&=\epsilon \gamma \sigma \left(1-Q^*-2I^*\right) Q^*-I^*C^*,
\end{align}
\end{subequations}
with initial conditions
\begin{equation}
 P^*(0)=1,\quad Q^*(0)=0,\quad C^*(0)=1, \quad  I^*(0)=\phi.\label{eq:nondim_init}
\end{equation}

The dimensionless parameters are then defined as,
\begin{equation}
\sigma=\frac{n_{0}}{c_0},\quad \rho=\frac{p_0}{c_0}, \quad \epsilon^2=\frac{k_2}{k_1},\quad \epsilon^2 \beta=\frac{k_4}{k_1}, \quad \epsilon \gamma=\frac{k_3}{k_1},\label{eq:dimpara}
\end{equation}
where \(\beta\), \(\gamma\), \(\sigma\) are order 1, and \(\epsilon\ll 1\) quantifies the disparity of the reaction rates. The magnitude of \(\rho\) will define \emph{moderate} (order 1) or \emph{high} (order \(\epsilon^{-1}\) )hydrogen peroxide regimes. Note that the definition of \(\epsilon\) is slightly different from that used in ref.~\cite{kerr2019mathematical}, which expressed the rates in terms of the rate constant \(k_0\) of a simplified form of the slow reaction. In what follows, asterisks will be dropped for brevity.

\begin{figure}
    \centering
    \includegraphics{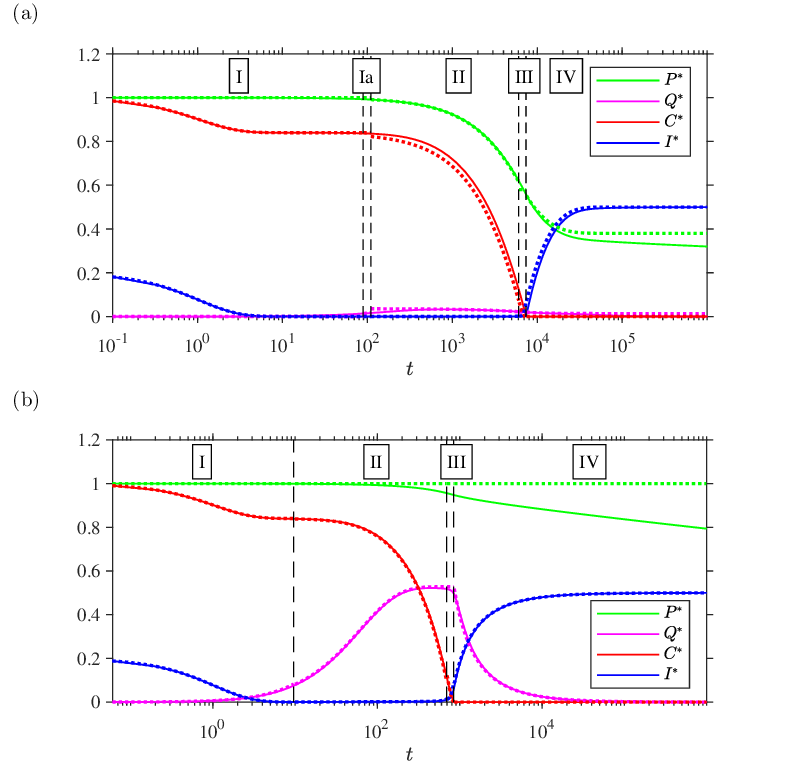}
    \caption{Comparison of numerical solutions (solid lines) and asymptotic solutions (dotted lines), as found in sections~\ref{sec:M-HP} and \ref{sec:H-HP} respectively to the dimensionless problem \eqref{eq:5fullmodelnond} plotted against logarithmic time. Parameter values (chosen arbitrarily): \(\epsilon=10^{-2}\), \(\beta=0.6\), \(\gamma=0.7\), \(\sigma=0.8\), \(\phi=0.2\). (a) Moderate hydrogen peroxide case M-HP with \(\rho=2\); (b) high hydrogen peroxide case H-HP with \(\hat{\rho}=0.9\) so that \(\rho=\hat{\rho}\epsilon^{-1}=90\).}
    \label{fig:NumericalSolutionLogRegions}
\end{figure}

As in earlier studies \cite{billingham1992mathematical,billingham1993mathematical,kerr2019mathematical}, we will exploit the small parameter \(\epsilon\) to develop an approximate solution to the problem via matched asymptotic expansions. The regimes of moderate and high hydrogen peroxide will be denoted \emph{M-HP} and \emph{H-HP}, and specifically refer to whether the parameter \(\rho\) is order 1 or order \(\epsilon^{-1}\), which determines whether the hypoiodous acid producing step has magnitude \(\epsilon^2\) or \(\epsilon\) in the dimensionless system.
We will now proceed to analyse each regime region-by-region. A comparison of numerical solutions to equations~\eqref{eq:5fullmodelnond} (computed with Matlab \texttt{ode23s}) and the asymptotic solutions which will be constructed is given in figure~\ref{fig:NumericalSolutionLogRegions}. 
The switchover time can be seen at the point where vitamin C (red line) drops, inducing a sharp increase in iodine (blue line).

\section{Model M-HP: moderate hydrogen peroxide \(\rho=O(1)\)}\label{sec:M-HP}
Under the assumption that \(\rho\) is order 1 (i.e.\ hydrogen peroxide concentration is initially comparable with vitamin C concentration), we note that hypoiodous acid \(Q\) is initially zero and its rate of production is order \(\epsilon\). Therefore it is appropriate to rescale \(Q=\epsilon\tilde{Q}\) where \(\tilde{Q}\) is order 1. Then the rescaled system is,
\begin{subequations}\label{eq:M-HPrscl}
\begin{align}
\frac{d P}{dt}&=-\epsilon^2 \sigma \left(1-\epsilon \tilde{Q}-2I\right)P-\epsilon^3 \beta \sigma \tilde{Q} P,\\
\frac{d \tilde{Q}}{dt}&=\epsilon \rho \left(1-\epsilon \tilde{Q}-2I\right)P- \epsilon \gamma \sigma \left(1-\epsilon \tilde{Q}-2I\right)\tilde{Q}-\epsilon^2 \beta \rho \tilde{Q} P, \label{eq:M-HP-b} \\
\frac{d C}{dt}&=-\sigma I C,\\
\frac{d I}{dt}&=\epsilon^2 \gamma \sigma \left(1-\epsilon \tilde{Q}-2I\right) \tilde{Q}-I C, \label{eq:M-HP-d}
\end{align}
\end{subequations}
with initial conditions \eqref{eq:nondim_init}.
The Model M-HP \eqref{eq:M-HPrscl} will be studied through four asymptotic regions in detail in the following sections.

\subsection{Region I: Initial Adjustment}\label{sec:M-HP-I}
During Region I, the reactants are mixed before the induction phase begins. It can be described as where the independent variable $t$ is order 1 ($t=O(1)$) and the dependent variables $D,P$ are order 1. Seeking a solution in the form of asymptotic expansions:
\begin{equation*}
P=P_{0}+\epsilon P_{1}+\ldots, \quad \tilde{Q}=\tilde{Q}_{0}+\epsilon \tilde{Q}_{1}+\ldots,\quad
C=C_{0}+\epsilon C_{1}+\ldots,\quad
\textrm{and}\quad I=I_{0}+\epsilon I_{1}+\ldots,
\end{equation*}
the leading order terms for the asymptotic expansion of the Model M-HP \eqref{eq:M-HPrscl} are:
\begin{subequations}
\begin{align}
\frac{dP_{0}}{dt}&=0,\label{eq:M-HP-I-a}\\
\frac{d\tilde{Q}_{0}}{dt}&=0,\label{eq:M-HP-I-b}\\
\frac{dC_{0}}{dt}&= -\sigma I_0 C_0,\label{eq:M-HP-I-c}\\
\frac{dI_{0}}{dt}&= - I_0 C_0,\label{eq:M-HP-I-d}
\end{align}
\end{subequations}
with the initial conditions:
\begin{equation}
P_0(0)=1,\quad \tilde{Q}_0(0)=0,\quad C_0(0)=1, \quad I_0(0)=\phi.
\end{equation}
Equations \eqref{eq:M-HP-I-a} and \eqref{eq:M-HP-I-b}, alongside their initial conditions yield
\begin{equation*}
P_0=1,\quad \tilde{Q}_0=0.
\end{equation*}
Dividing \eqref{eq:M-HP-I-c} by \eqref{eq:M-HP-I-d} gives
\begin{equation}
    \frac{dC_0}{dI_0}=\sigma,
\end{equation}
which has solution
\begin{equation}
C_0=\sigma \left(I_0-\phi\right)+1\label{eq:M-HP-C0}.
\end{equation}
Substituting into equation~\eqref{eq:M-HP-I-d} then gives the differential equation in terms of \(I_0\) only,
\begin{equation*}
    \frac{dI_0}{dt}=-I_0\left[\sigma \left(I_0-\phi\right)+1\right],
\end{equation*}
which has solution,
\begin{equation}
    I_0=\frac{\phi \left(1-\sigma \phi\right) e^{-\left(1-\sigma \phi\right)t}}{1-\sigma \phi e^{-\left(1-\sigma \phi\right)t}}.\label{eq:M-HP-I0}
\end{equation}
Substituting \eqref{eq:M-HP-I0} into \eqref{eq:M-HP-C0} leads to
\begin{equation}
    C_0=\frac{\sigma \phi\left(1-\sigma \phi\right) e^{-\left(1-\sigma \phi\right)t}}{1-\sigma \phi e^{-\left(1-\sigma \phi\right)t}}+\left(1-\sigma \phi\right).
\end{equation}
Throughout it will be assumed that in dimensional variables, \(\iota_0<c_0\), which corresponds to the nondimensional condition \(\sigma \phi<1\). This condition corresponds to the initial vitamin C concentration being sufficient to survive the initial adjustment. The initial adjustment therefore corresponds to consumption of a proportion \(\sigma\phi\) of the initial vitamin C concentration via the fast reaction, thereby converting the bulk of the initial iodine concentration to iodide in order 1 time; the slow reaction is subleading in this region. In particular, note that the order 1 part of \(I_0(t)\) is exponentially decaying.

\subsection{Region Ia: Quasi Equilibrium Established for \(Q\)}
Following the exponential decay of vitamin C (to \(1-\sigma\phi\)) on the previous region, the fast reaction no longer dominates and the slow reaction becomes of comparable importance. By inspecting equation~\eqref{eq:M-HP-d} it is apparent that balance between fast and slow reactions implies that \(I(t)=O(\epsilon^2)\).

A distinctive feature of the M-HP case as opposed to the H-HP case is the existence of a region intermediate between the initial adjustment and the induction period, during which the hypoiodous acid concentration reaches quasi-equilibrium following the rise in iodide concentration. For consistency with ref.\ \cite{kerr2019mathematical} which used the designation `region II' specifically for the induction period, we will designate the intermediate region between the initial adjustment and the induction period by `Ia'. Physically, region Ia corresponds to the two forward steps in the slow reaction evolving towards a balance in their production and removal of hypoiodous acid. Inspecting equation~\eqref{eq:M-HP-b} shows that \(\tilde{Q}(t)\) varies on a timescale \(t=O(\epsilon^{-1})\). Therefore the scalings for region Ia are, $P=C=O(1)$, $I=O(\epsilon^2)$ and $t=O(\epsilon^{-1})$. Denoting $\tau=\epsilon t$ and $I=\epsilon^2 \tilde{I}$, then the system takes the form
\begin{subequations}
    \begin{align}
        \frac{d P}{d\tau}&=-\epsilon \sigma \left(1-\epsilon \tilde{Q}-2\epsilon^2 \tilde{I}\right)P-\epsilon^2 \beta \sigma \tilde{Q} P,\\
        \frac{d \tilde{Q}}{d\tau}&= \rho \left(1-\epsilon \tilde{Q}-2\epsilon^2 \tilde{I}\right)P - \gamma \sigma \left(1-\epsilon \tilde{Q}-2\epsilon^2 \tilde{I}\right)\tilde{Q}-\epsilon \beta \rho \tilde{Q} P,\\
        \frac{d C}{d\tau}&=-\epsilon\sigma \tilde{I} C,\label{eq:M-MP-Ia-c}\\
        \epsilon\frac{d \tilde{I}}{d\tau}&= \gamma \sigma \left(1-\epsilon \tilde{Q}-2\epsilon^2 \tilde{I}\right) \tilde{Q} -  \tilde{I} C. \label{eq:M-HP-Ia-d}
    \end{align}
\end{subequations}
The leading order terms of the asymptotic expansions therefore satisfy,
\begin{subequations}\label{eq:M-HP-Ia}
\begin{align}
\frac{dP_{0}}{d\tau}&=0,\label{eq:M-HP-Ia-leading-a}\\
\frac{d\tilde{Q}_{0}}{d\tau}&=\rho P_0-\gamma \sigma \tilde{Q}_0,\label{eq:M-HP-Ia-leading-b}\\
\frac{dC_{0}}{d\tau}&= 0, \label{eq:M-HP-Ia-leading-c}\\
0&= \gamma \sigma \tilde{Q}_0-\tilde{I}_0 C_0,\label{eq:M-HP-Ia-leading-d}
\end{align}
\end{subequations}
with the matching conditions to the Region I solution:
\begin{equation*}
P_0(0)=1,\quad \tilde{Q}_0(0)=0,\quad C_0(0)=1-\sigma \phi, \quad  \tilde{I}_0(0)=0.
\end{equation*}
Solving equations~\eqref{eq:M-HP-Ia-leading-a} and \eqref{eq:M-HP-Ia-leading-c} then yields the constant solutions,
\begin{equation}\label{eq:M-HP-Ia-PCsol}
    P_0(\tau)=1,\quad C_0(\tau)=1-\sigma \phi.
\end{equation}
Equations~\eqref{eq:M-HP-Ia-leading-b} and \eqref{eq:M-HP-Ia-leading-d} can then be solved to yield,
\begin{equation}
\tilde{Q}_0=\frac{\rho}{\gamma \sigma}\left(1-e^{-\gamma \sigma \tau}\right),\label{eq:M-HP-Ia-Q0sol}
\end{equation}
and
\begin{equation}
\tilde{I}_0=\frac{\rho}{1-\sigma \phi}\left(1-e^{-\gamma \sigma \tau}\right).
\end{equation}
The limiting behaviours of the solution in region Ia are therefore,
\begin{equation}
    P_0(\tau) = 1, \quad \tilde{Q}_0(\tau) \rightarrow \frac{\rho}{\gamma\sigma}, \quad C_0(\tau) = 1-\sigma\phi, \quad \tilde{I}_0(\tau) \rightarrow \frac{\rho}{1-\sigma\phi},
\end{equation}
as \(\tau\rightarrow \infty\).

\subsection{Region II: induction period}
Following the establishment of quasi-equilibrium within the slow reaction, the system enters the induction period, which is the distinguishing feature of a clock reaction. During the induction period, the fast and slow reactions are in quasi-equilibrium, with iodine levels held low due to the disparity  in reaction rates. The vitamin C is gradually consumed until it is exhausted and the switchover process then occurs. From inspecting equation~\eqref{eq:M-MP-Ia-c} it is apparent that this process occurs on a timescale of order \(\epsilon^{-1}\) slower than \(\tau\), motivating the introduction of a second rescaled time variable \(T=\epsilon \tau = \epsilon^{2} t\). The rescaled system is then,
\begin{subequations}
    \begin{align}
        \frac{d P}{dT}&=- \sigma \left(1-\epsilon \tilde{Q}-2\epsilon^2 \tilde{I}\right)P - \epsilon \beta \sigma \tilde{Q} P,\\
        \epsilon\frac{d \tilde{Q}}{dT}&= \rho \left(1-\epsilon \tilde{Q}-2\epsilon^2 \tilde{I}\right)P - \gamma \sigma \left(1-\epsilon \tilde{Q}-2\epsilon^2 \tilde{I}\right)\tilde{Q}-\epsilon \beta \rho \tilde{Q} P,\\
        \frac{d C}{dT}&=-\sigma \tilde{I} C,\\
        \epsilon^2 \frac{d \tilde{I}}{dT}&= \gamma \sigma \left(1-\epsilon \tilde{Q}-2\epsilon^2 \tilde{I}\right) \tilde{Q} - \tilde{I} C.
    \end{align}
\end{subequations}
The leading order terms of the asymptotic expansions in region II therefore satisfy,
\begin{subequations}\label{eq:M-HP-II}
\begin{align}
\frac{dP_{0}}{dT}&=-\sigma P_{0},\label{eq:M-HP-II-a}\\
0&=\rho P_0-\gamma \sigma \tilde{Q}_0,\label{eq:M-HP-II-b}\\
\frac{dC_{0}}{dT}&= -\sigma \tilde{I}_0C_0,\label{eq:M-HP-II-c}\\
0&= \gamma \sigma \tilde{Q}_0-\tilde{I}_0 C_0.\label{eq:M-HP-II-d}
\end{align}
\end{subequations}
Solving equations~\eqref{eq:M-HP-II} and matching to the region Ia solutions then yields the leading order approximations,
\begin{subequations}
    \begin{align}
        P_0(T)         & = e^{-\sigma T}, \\
        \tilde{Q}_0(T) & = \frac{\rho}{\gamma \sigma}e^{-\sigma T}, \\
        C_0(T)         & = \rho e^{-\sigma T}+(1-\sigma \phi-\rho), \label{eq:M-HP-II_summary-c}\\
        \tilde{I}_0(T) & = \frac{\rho e^{-\sigma T}}{\rho e^{-\sigma T}+(1-\sigma \phi-\rho)}.
    \end{align}
\end{subequations}


The switchover time can then be approximated as the time at which the region II solution predicts that the vitamin C concentration is zero at leading order. Inspecting equation~\eqref{eq:M-HP-II_summary-c} yields in dimensionless variables,
\begin{equation}
    T_{sw} = \frac{1}{\sigma}\ln{\left(\frac{\rho}{\rho+\sigma\phi-1}\right)}, \label{eq:M-HP-Tsw-nondim}
\end{equation}
or in dimensional variables,
\begin{equation}
t_{sw}^{\text{M-HP}}=\frac{1}{k_2 n_0}\ln{\left(\frac{p_0}{p_0+\phi n_0-c_0}\right)}. \label{eq:M-HP-tsw}
\end{equation}
This equation emphasises that there must be sufficient combined hydrogen peroxide and molecular iodine provided initially to drive the vitamin C concentration to zero before the hydrogen peroxide is exhausted, i.e.\ \(p_0 + \phi n_0 > c_0\).

A similar equation to \eqref{eq:M-HP-tsw} was given in ref.\ \cite{parra2020enhancing} through a quasi-steady argument and making the phenomenological assumption that the slow reaction is linear in iodide concentration, (although not including the \(\phi n_0\) term in the denominator). However, by accounting for the sub-steps of the slow reaction, the dynamics of this regime emerge naturally from the law of mass action. Equation~\eqref{eq:M-HP-tsw} is one of the two formulae that will be tested experimentally in this paper.

While the above analysis provides a formula for the switchover time, as in ref.\ \cite{kerr2019mathematical} it is of interest to construct approximate solutions for the remaining dynamics of the system. 

\subsection{Region III: Corner}\label{sec:M-HP-III}
As the vitamin C concentration is exhausted, the iodine concentration begins to change significantly, and the system is no longer in equilibrium. Following ref.~\cite{kerr2019mathematical} we will refer to this as the `corner' region; the analysis is similar, with some minor adaptations for the four variable model and different definition of the small parameter \(\epsilon\).

Shifting the time coordinate so that the origin corresponds to the switchover time, the most structured balance in equation~\eqref{eq:M-HPrscl} is given by the `crossover' scalings,
\begin{equation}
    \bar{t}=\epsilon^{-1}(\epsilon^2 t - T_{sw}), \quad \bar{C}=\epsilon^{-1}C, \quad \bar{I}=\epsilon^{-1} I, \quad \bar{Q}=\tilde{Q} = \epsilon^{-1}Q.
\end{equation}
The rescaled system is then,
\begin{subequations}
    \begin{align}
        \frac{d P}{d\bar{t}}&=- \epsilon\sigma \left(1-\epsilon \tilde{Q}-2\epsilon \bar{I}\right)P - \epsilon^2 \beta \sigma \tilde{Q} P,\\
        \frac{d \tilde{Q}}{d\bar{t}}&= \rho \left(1-\epsilon \tilde{Q}-2\epsilon\bar{I}\right)P - \gamma \sigma \left(1-\epsilon \tilde{Q}-2\epsilon \bar{I}\right)\tilde{Q}-\epsilon \beta \rho \tilde{Q} P,\\
        \frac{d \bar{C}}{d\bar{t}}&=-\sigma \bar{I} \bar{C},\\
        \frac{d \bar{I}}{d\bar{t}}&= \gamma \sigma \left(1-\epsilon \tilde{Q}-2\epsilon\bar{I}\right) \tilde{Q} - \bar{I} \bar{C}.
    \end{align}
\end{subequations}

Taking the leading order terms of the asymptotic expansions yields,
 \begin{subequations}\label{eq:M-HP-III-leading}
\begin{align}
\frac{d P_{0}}{d\bar{t}}&=0,\label{eq:M-HP-III-leading-a}\\
\frac{d\tilde{Q}_{0}}{d\bar{t}}&=\rho P_0-\gamma \sigma \tilde{Q}_0,\label{eq:M-HP-III-leading-b}\\
\frac{d\bar{C}_{0}}{d\bar{t}}&= -\sigma \bar{I}_0\bar{C}_0,\label{eq:M-HP-III-leading-c}\\
\frac{d\bar{I}_{0}}{d\bar{t}}&= \gamma \sigma \tilde{Q}_0-\bar{I}_0\bar{C}_0.\label{eq:M-HP-III-leading-d}
\end{align}
\end{subequations}

Matching~\eqref{eq:M-HP-III-leading-a} to region II at \(T=T_{sw}\) gives \(P_0(\bar{t})=(\rho+\sigma\phi-1)/\rho\). Substituting into equation~\eqref{eq:M-HP-III-leading-b} and solving for \(\bar{Q}_0\), we arrive at
\begin{equation}
    \bar{Q}_0(\bar{t}) = \frac{\rho+\sigma\phi-1}{\gamma\sigma} + c_1 e^{-\gamma\sigma \bar{t}}, \label{eq:H-HP-III-Q0sol}
\end{equation}
where \(c_1\) is a constant. Matching to region II as \(\bar{t}\rightarrow -\infty\) establishes that \(c_1=0\) so that \(\bar{Q}_0\) remains finite; indeed the resulting constant solution \(\bar{Q}_0(\bar{t})=(\rho+\sigma\phi-1)/(\gamma\sigma)\) then matches with the region II solution at \(T=T_{sw}\).

Subtracting a multiple of equation~\eqref{eq:M-HP-III-leading-d} from \eqref{eq:M-HP-III-leading-c} and substituting the solution for \(\bar{Q}_0\) then gives the equation,
\begin{equation}
    \frac{d}{d\bar{t}}(\bar{C}_0 - \sigma \bar{I}_0) = -\sigma(\rho+\sigma\phi-1),
\end{equation}
from which it follows that
\begin{equation}
    \bar{C}_0(\bar{t}) - \sigma\bar{I}_0(\bar{t}) = -\sigma(\rho+\sigma\phi-1)\bar{t} + c_2, \label{eq:M-HP-III-cons}
\end{equation}
where \(c_2\) is a constant. Noting that an order 1 change to \(c_2\) corresponds to an order 1 shift in the origin for \(\bar{t}\), and hence a subleading change to \(t_{sw}\), it is reasonable to choose \(c_2=0\) for convenience in the present analysis. A more precise choice of \(c_2\) would require the next order solution for \(C(t)\) in region II to be calculated for the purpose of matching. Using equation~\eqref{eq:M-HP-III-cons} to substitute for \(\bar{I}_0\) in equation~\eqref{eq:M-HP-III-leading-c} then yields a Riccati equation for \(\bar{C}_0(\bar{t})\),
\begin{equation}
    \frac{d\bar{C}_{0}}{d\bar{t}}= -\left(\sigma\left(\rho+\sigma \phi-1\right) \bar{t}+\bar{C}_0\right)\bar{C}_0,
\end{equation}
as found for example in the study of a cubic autocatalysis clock reaction \cite{billingham1992mathematical}. It is convenient to introduce the notation \(\bar{\rho}=\sqrt{\sigma(\rho+\sigma\phi-1)}\) which renders the equation into an identical form to that in ref.~\cite{kerr2019mathematical}, from which we can immediately write down the solution,
\begin{equation}
    \bar{C}_0=\frac{\bar{\rho}\sqrt{2}\exp{[-\bar{\rho}^2\bar{t}^2/2 ] }}{\sqrt{\pi}\left[\mathrm{erf}\left(\bar{\rho}\bar{t}/\sqrt{2}\right)+c_4\right]} ,\label{eq:M-HP-III-C0sol}
\end{equation}
where \(c_4\) is another constant to be determined. The constant \(c_4\) can be deduced by examining the limiting behaviour of the error function as \(\bar{t}\rightarrow -\infty\); for the solution for \(\bar{C}_0\) to approach a non-zero limit requires that \(c_4=1\). The solution for \(\bar{I}_0(\bar{t})\) follows immediately from equation~\eqref{eq:M-HP-III-cons}. In summary, the solutions for region III are,
\begin{subequations}
    \begin{align}
        \bar{P}_0(\bar{t}) & = \frac{\bar{\rho}^2}{\rho\sigma}, \\
        \bar{Q}_0(\bar{t}) & = \frac{\bar{\rho}^2}{\gamma\sigma^2}, \\
        \bar{C}_0(\bar{t}) & = \frac{\bar{\rho}\sqrt{2}\exp{[-\bar{\rho}^2\bar{t}^2/2 ] }}{\sqrt{\pi}\left[\mathrm{erf}\left(\bar{\rho}\bar{t}/\sqrt{2}\right)+1\right]}, \\
        \bar{I}_0(\bar{t}) & = \frac{\bar{\rho}\sqrt{2}\exp{[-\bar{\rho}^2\bar{t}^2/2 ] }}{\sqrt{\pi}\sigma\left[\mathrm{erf}\left(\bar{\rho}\bar{t}/\sqrt{2}\right)+1\right]} + \frac{\bar{\rho}^2}{\sigma}\bar{t}.
    \end{align}
\end{subequations}

\subsection{Region IV: equilibration}
Following the exhaustion of vitamin C and the switchover, the system approaches its long term equilibrium state. 
From equations~\eqref{eq:5fullmodelnond} it is clear that equilibrium with \(P>0\) and \(C=0\) also requires that \(Q=0\) and \(I=1/2\) -- in other words, all iodide and hypoiodous acid is ultimately converted to molecular iodine.

Using the same timescale \(T=\epsilon^2 t\) as the induction period, but with \(I(T)\) now order 1, the rescaled system is,
\begin{subequations}
    \begin{align}
        \epsilon^2 \frac{dP}{dT} & = -\epsilon^2 \sigma (1-\epsilon\tilde{Q}-2I)P - \epsilon^3 \beta \sigma \tilde{Q}P, \label{eq:M-HP-IV-a}\\
        \epsilon^2 \frac{d\tilde{Q}}{dT} & = \epsilon \rho(1-\epsilon\tilde{Q}-2I) P - \epsilon\gamma\sigma(1-\epsilon\tilde{Q}-2I)\tilde{Q}-\epsilon^2\beta\rho\tilde{Q}P, \label{eq:M-HP-IV-b}\\
        \epsilon^2 \frac{dC}{dT} & = -\sigma I C, \label{eq:M-HP-IV-c}\\
        \epsilon^2 \frac{dI}{dT} & = \epsilon^2 \gamma \sigma(1-\epsilon\tilde{Q}-2I)\tilde{Q} - I C.\label{eq:M-HP-IV-d}
    \end{align}
\end{subequations}
Equation~\eqref{eq:M-HP-IV-c} yields that \(I_0 C_0=0\) at leading order, from which we deduce the expected result that \(C_0(T)=0\) for all \(T\). Indeed it can be shown by mathematical induction that \(C_n(T)=0\) at all powers of \(n\), by a similar argument to that given by Kerr et al.~\cite{kerr2019mathematical}. The leading order terms of the rest of the system are,
\begin{subequations}
    \begin{align}
        \frac{dP_0}{dT} & = -\sigma (1-2I_0)P_0, \label{eq:M-HP-IV-P0}\\
        0 = & \rho (1-2I_0)P_0 - \gamma\sigma(1-2I_0)\tilde{Q}_0, \\
        \frac{dI_0}{dt} & = \gamma\sigma(1-2I_0)\tilde{Q}_0.
    \end{align}
\end{subequations}
Rearranging this system yields,
\begin{equation}
    \frac{dP_0}{dI_0} = -\frac{\sigma}{\rho}, \quad \mbox{hence} \quad P_0(T) = -\frac{\sigma}{\rho} I_0 + b_1, \label{eq:M-HP-IV-P-I}
\end{equation}
where \(b_1\) is a constant. Matching to region II immediately yields \(b_1 = \exp(-\sigma T_{sw})=(\rho+\sigma\phi-1)/\rho\), i.e.\ the level of hydrogen peroxide remaining after the induction period. We will retain the notation \(b_1\) for the rest of the calculation for brevity. Substituting for \(I_0\) in equation~\eqref{eq:M-HP-IV-P0} then yields an equation involving only \(P_0\),
\begin{equation}
    \frac{dP_0}{dT} = P_0\left(-2\rho P_0 + 2\rho b_1 - \sigma\right),
\end{equation}
which integrates to,
\begin{equation}
    \ln\left|\frac{2\rho P_0 - (2\rho b_1 -\sigma)}{P_0}\right| = -(2\rho b_1 -\sigma)(T+b_2),    
\end{equation}
where \(b_2\) is another constant to be determined by matching. Noting that \(P_0\geqslant 0\) and moreover at the beginning of region IV, \(P_0=b_1\), it follows that the term inside the modulus signs is initially positive; a change of sign of the numerator cannot then occur at finite time, so we deduce that,
\begin{equation}
    \ln\left(\frac{2\rho P_0 - (2\rho b_1 -\sigma)}{P_0}\right) = -(2\rho b_1 -\sigma)(T+b_2). \label{eq:M-HP-IV-lnP}
\end{equation}
The constant \(b_2\) must then be,
\begin{equation}
    b_2=-T_{sw}+\frac{1}{2\rho b_1-\sigma} \ln\left(\frac{b_1}{\sigma}\right).
\end{equation}
Rearranging equation~\eqref{eq:M-HP-IV-lnP} then gives the long-timescale solution for hydrogen peroxide,
\begin{equation}
    P_0(T) = \frac{2\rho b_1-\sigma}{2\rho-(\sigma/b_1)\exp(-(2\rho b_1-\sigma)(T-T_{sw}))},\quad \mbox{where} \quad b_1=\frac{\rho+\sigma\phi-1}{\rho}, \label{eq:M-HP-IV-Psol}
\end{equation}
and substitution into equation~\eqref{eq:M-HP-IV-P-I} shows moreover that,
\begin{equation}
    I_0(T) = \rho b_1\left(\frac{1-\exp(-(2\rho b_1-\sigma)(T-T_{sw}))}{2\rho b_1-\sigma\exp(-(2\rho b_1-\sigma)(T-T_{sw}))}\right).
\end{equation}

Two cases can be distinguished in equation~\eqref{eq:M-HP-IV-Psol}. If \(2\rho b_1-\sigma>0\), i.e.\ there is sufficient hydrogen peroxide remaining following the switchover, then in the limit as \(T\rightarrow \infty\),
\begin{subequations}
    \begin{align}
        P_0(T) & \sim \frac{2\rho b_1-\sigma}{2\rho}=\frac{2(\rho-1)+\sigma(2\phi-1)}{2\rho}, \\
        I_0(T) & \sim \frac{1}{2},
    \end{align}
\end{subequations}
in other words, all of the iodide and hypoiodous acid is converted to molecular iodine.

If however \(2\rho b_1 -\sigma<0\), i.e.\ there is insufficient hydrogen peroxide remaining following the switchover, then as \(T\rightarrow \infty\),
\begin{subequations}
    \begin{align}
        P_0(T) & \sim 0, \\
        I_0(T) & \sim \frac{\rho b_1}{\sigma}=\frac{\rho+\sigma\phi-1}{\sigma}.
    \end{align}
\end{subequations}
In dimensional variables, the condition \(2\rho b_1-\sigma> 0\) is equivalent to \(p_0> n_0/2 + c_0 - \iota_0\). To convert all iodide to iodine by the end of the process, there must be sufficient hydrogen peroxide to account for all supplied atomic iodine (noting one hydrogen peroxide molecule is required for every two iodide ions), plus exhausting all supplied vitamin C, minus the initial adjustment.

In summary, the region IV leading order solutions describing the long-term evolution of the system are,
\begin{subequations}
    \begin{align}
    P_0(T) & = \frac{2\rho b_1-\sigma}{2\rho-(\sigma/b_1)\exp(-(2\rho b_1-\sigma)(T-T_{sw}))}, \\
    \tilde{Q}_0(T) & = \frac{\rho}{\gamma\sigma}\left[\frac{2\rho b_1-\sigma}{2\rho-(\sigma/b_1)\exp(-(2\rho b_1-\sigma)(T-T_{sw}))}\right],\\
    C_0(T) & = 0, \\
    I_0(T) & = \rho b_1\left[\frac{1-\exp(-(2\rho b_1-\sigma)(T-T_{sw}))}{2\rho b_1-\sigma\exp(-(2\rho b_1-\sigma)(T-T_{sw}))}\right],
    \end{align}
\end{subequations}
where \(T_{sw}\) is given in equation~\eqref{eq:M-HP-Tsw-nondim} and \(b_1=(\rho+\sigma\phi-1)/\rho\).

\begin{figure}
    \centering
    \includegraphics{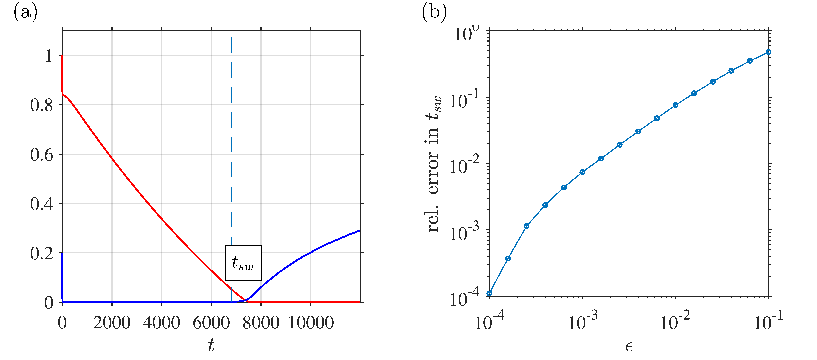}
    \caption{Comparison of asymptotic-derived formula and numerical solution in determination of the dimensionless switchover time for regime M-HP. 
    Parameter values (chosen arbitrarily): \(\epsilon=10^{-2}\), \(\beta=0.6\), \(\gamma=0.7\), \(\sigma=0.8\), \(\phi=0.2\) and \(\rho=2\). (a) Numerical solutions for \(C(t)\) and \(I(t)\) shown alongside the value of the switchover time (dashed line) calculated by the formula \eqref{eq:M-HP-Tsw-nondim}. (b) Relative error between the asymptotic expression  \eqref{eq:M-HP-Tsw-nondim} and the point at which the numerical solution falls below a threshold value of \(\epsilon\), as a function of the small parameter \(\epsilon\). }\label{fig:M-HP-Convergence}
\end{figure}

A comparison between numerical and leading order asymptotic solutions in each region in dimensionless variables is given in figure~\ref{fig:NumericalSolutionLogRegions}(a), for a specific set of parameters, with small parameter chosen as \(\epsilon=10^{-2}\). In the absence of detailed measurement of the reaction rates in the system, this choice is arbitrary. Additional detail on the dynamics for \(Q(t)\) in region Ia is given in Appendix~\ref{sec:appendix}, figure~\ref{fig:Comparison-M-Ia}. As shown in figure~\ref{fig:NumericalSolutionLogRegions}(a), there is a small but noticeable difference between the asymptotic and numerical solutions as \(C(t)\) approaches zero in the induction region which results in a small discrepancy in switchover time (figure~\ref{fig:M-HP-Convergence}(a)). By plotting as a function of \(\epsilon\) the relative difference in the asymptotic switchover time value (equation~\eqref{eq:M-HP-Tsw-nondim}) and the point at which the numerical solution falls below a threshold value of \(\epsilon^{-1}\) (figure~\ref{fig:M-HP-Convergence}(b)), this discrepancy converges to zero. As would be expected from taking only the leading order terms of the asymptotic expansion, the convergence is approximately linear in \(\epsilon\).

\section{Model H-HP: high hydrogen peroxide \(\rho=O(\epsilon^{-1})\)}\label{sec:H-HP}
The contrasting regime with high hydrogen peroxide can be characterised mathematically by the ratio \(p_0/c_0 = \rho \gg 1\), or more precisely, \(\rho\) being order \(\epsilon^{-1}\). Therefore we will write \(\rho = \epsilon^{-1} \hat{\rho}\), where \(\hat{\rho}\) is order 1. The dimensionless model of equation~\eqref{eq:5fullmodelnond} then becomes,
\begin{subequations}\label{eq:H-HP}
\begin{align}
\frac{d P}{dt}&=-\epsilon^2 \sigma \left(1- Q-2I\right)P-\epsilon^2 \beta \sigma Q P,\\
\frac{d Q}{dt}&=\epsilon \hat{\rho} \left(1- Q-2I\right)P- \epsilon \gamma \sigma \left(1- Q-2I\right)Q-\epsilon \beta \hat{\rho} Q P,\\
\frac{d C}{dt}&=-\sigma I C,\label{eq:H-HP-c}\\
\frac{d I}{dt}&=\epsilon \gamma \sigma \left(1-Q-2I\right) Q-I C, \label{eq:H-HP-d}
\end{align}
\end{subequations}
with the same initial conditions \eqref{eq:nondim_init} as previously. 

\subsection{Region I: Initial Adjustment}\label{sec:H-HP-I}
The initial process of vitamin C being rapidly consumed to convert the bulk of the iodide to molecular iodine occurs during a time interval of order 1 similarly to the M-HP regime.
Again it is necessary to rescale the hypoiodous acid concentration, so that \(Q=\epsilon\tilde{Q}\), with \(\tilde{Q}\) being order 1. Then the rescaled system is,
\begin{subequations}\label{eq:H-HP-I}
\begin{align}
\frac{d P}{dt}&=-\epsilon^2 \sigma \left(1- \epsilon\tilde{Q}-2I\right)P-\epsilon^3 \beta \sigma \tilde{Q} P,\\
\epsilon\frac{d \tilde{Q}}{dt}&=\epsilon \hat{\rho} \left(1- \epsilon\tilde{Q}-2I\right)P- \epsilon^2 \gamma \sigma \left(1- \epsilon\tilde{Q}-2I\right)\tilde{Q}-\epsilon^2 \beta \hat{\rho} \tilde{Q} P,\\
\frac{d C}{dt}&=-\sigma I C,\\
\frac{d I}{dt}&=\epsilon^2 \gamma \sigma \left(1-\epsilon\tilde{Q}-2I\right) \tilde{Q}-I C.
\end{align}
\end{subequations}

The leading order system becomes,
\begin{subequations}
\begin{align}
\frac{dP_{0}}{dt}&=0,\\
\frac{d\tilde{Q}_{0}}{dt}&=\hat{\rho}\left(1-2I_0\right)P_0,\label{eq:H-HP-I-b}\\
\frac{dC_{0}}{dt}&= -\sigma I_0 C_0,\label{eq:H-HP-I-c}\\
\frac{dI_{0}}{dt}&= - I_0 C_0.\label{eq:H-HP-I-d}
\end{align}
\end{subequations}
The solutions for \(P_0\), \(C_0\) and \(I_0\) are identical to those in section~\ref{sec:M-HP}\ref{sec:M-HP-I}. The problem for \(\tilde{Q}_0\) then reduces to the separable equation,
\begin{equation}
   \frac{d\tilde{Q}_0}{dt}=\hat{\rho}\bigg[1-\frac{2\phi \left(1-\sigma \phi\right) e^{-\left(1-\sigma \phi\right)t}}{1-\sigma \phi e^{-\left(1-\sigma \phi\right)t}}\bigg], 
\end{equation}
with solution,
\begin{equation}
    \tilde{Q}_0=\hat{\rho}t-\frac{2\hat{\rho}}{\sigma}\ln{\left(1-\sigma \phi e^{-\left(1-\sigma \phi\right)t}\right)}+\frac{2\hat{\rho}}{\sigma}\ln{\left(1-\sigma \phi\right)}.
\end{equation}
The limiting behaviour as \(t\rightarrow \infty\) is then,
\begin{equation}
    P_0(t)\rightarrow 1, \quad \tilde{Q}_0(t)\sim \hat{\rho}t, \quad C_0(t)\rightarrow 1-\sigma \phi, \quad I_0 \rightarrow 0.
\end{equation}
The linear growth of \(\tilde{Q}_0(t)\) implies that \(Q(t)\) becomes order 1 on a timescale of order \(\epsilon^{-1}\).

The more rapid growth of \(Q(t)\), arising because the hypoiodous acid producing step is no longer rate limiting in case of abundant hydrogen peroxide, means that it is no longer necessary to consider a separate region for quasi-equilibrium to be established, so there is no counterpart to region Ia studied for the M-HP regime. 

\subsection{Region II: Induction period}\label{sec:H-HP-II}
Moving on to region II -- the induction period -- it is apparent that quasi-equilibrium for forward and backward reactions occurs in equation~\eqref{eq:H-HP-d} when \(I(t)\) is order \(\epsilon\), i.e.\ an order of magnitude larger than in the moderate hydrogen peroxide case. Considering equation~\eqref{eq:H-HP-c} with the rescaling \(I(t)=\epsilon\tilde{I}(t)\) (with \(\tilde{I}(t)\) being order 1), it is clear that the appropriate time rescaling is then \(t=\epsilon^{-1}\tau\), with \(\tau\) being order 1. The resulting system for region II is then,
\begin{subequations}\label{eq:H-HP-II}
\begin{align}
\frac{d P}{d\tau}&=-\epsilon \sigma \left(1- Q-2\epsilon\tilde{I}\right)P-\epsilon \beta \sigma Q P,\\
\frac{d Q}{d\tau}&= \hat{\rho} \left(1- Q-2\epsilon\tilde{I}\right)P- \gamma \sigma \left(1- Q-2\epsilon\tilde{I}\right)Q- \beta \hat{\rho} Q P,\\
\frac{d C}{d\tau}&=-\sigma \tilde{I} C,\label{eq:H-HP-II-c}\\
\epsilon\frac{d \tilde{I}}{d\tau}&= \gamma \sigma \left(1-Q-2\epsilon\tilde{I}\right) Q-\tilde{I} C. \label{eq:H-HP-II-d}
\end{align}
\end{subequations}
The leading order system takes the form,
\begin{subequations}\label{eq:H-HP-II-leading}
\begin{align}
\frac{dP_{0}}{d\tau}&=0,\label{eq:H-HP-II-leading-a}\\
\frac{dQ_{0}}{d\tau}&=\hat{\rho} \left(1-Q_0\right)P_0-\gamma \sigma \left(1-Q_0\right)Q_0-\beta\hat{\rho}Q_0P_0,\label{eq:H-HP-II-leading-b}\\
\frac{dC_{0}}{d\tau}&= -\sigma \tilde{I}_0C_0,\label{eq:H-HP-II-leading-c}\\
0&= \gamma \sigma \left(1-Q_0\right)Q_0-\tilde{I}_0 C_0.\label{eq:H-HP-II-leading-d}
\end{align}
\end{subequations}

After integrating equation~\eqref{eq:H-HP-II-leading-a} and matching to region I, the hydrogen peroxide concentration remains constant at leading order, i.e.\ \(P_0(\tau)=1\). Substituting into equation~\eqref{eq:H-HP-II-leading-b} and solving yields,
\begin{equation}
    Q_0(\tau)=-\frac{b}{\gamma \sigma}\left[\frac{\tanh\left(b\tau\right)+\tanh\left(bc_5\right)}{1+\tanh\left(b\tau\right)\tanh\left(bc_5\right)}\right]+\frac{a}{2\gamma \sigma}.
\end{equation}
where the constants \(a\) and \(b\) are given by,
\begin{equation}
a=\hat{\rho}+\gamma\sigma+\beta \hat{\rho}\quad\text{and}\quad 
b=\frac{\sqrt{-4\gamma\sigma\hat{\rho}+a^2}}{2},\label{eq:H-HP-II-const-b}
\end{equation}
and \(c_5\) is to be determined via matching. Because \(a^2-4\gamma\sigma\hat{\rho}=(\hat{\rho}(\beta-1)+\gamma\sigma)^2+4\beta\hat{\rho}^2>0\), it can be deduced that \(b\) is real-valued.

Matching to region I (as $t\rightarrow\infty$, $Q_0(t)\rightarrow 0$), gives
\begin{equation}
    c_5=\frac{1}{b}\mathrm{arctanh}\left(\frac{a}{2b}\right),\label{eq:H-HP-II-c1}
\end{equation}
from which we deduce the solution,
\begin{equation}
    Q_0(\tau)=-\frac{b}{\gamma \sigma}\left[\frac{2b\tanh\left(b\tau\right)+a}{2b+a\tanh\left(b\tau\right)}\right]+\frac{a}{2\gamma \sigma}.\label{eq:H-HP-II-Q0}
\end{equation}

Combining equations~\eqref{eq:H-HP-II-leading-c}, \eqref{eq:H-HP-II-leading-d} and substituting the solution for \(Q_0\) then reduces the problem for \(C_0\) to
\begin{equation}
    \frac{dC_0}{d\tau}=\sigma \left(1+\frac{b}{\gamma \sigma}
    \left[\frac{2b\tanh\left(b\tau\right)+a}{2b+a\tanh\left(b\tau\right)}\right]
    -\frac{a}{2\gamma \sigma}\right)
    \left(b\left[\frac{2b\tanh\left(b\tau\right)+a}{2b+a\tanh\left(b\tau\right)}\right]-\frac{a}{2}\right)
    , \label{eq:H-HP-II-dC0}
\end{equation}
which integrates to yield,
\begin{equation}
    C_0(\tau)=w_1\ln{\bigg\lvert \cosh(b\tau)+\frac{a}{2b}\sinh(b\tau)\bigg\rvert}+w_{2}\frac{2b}{2b+a\tanh\left(b\tau\right)}+w_{3}\tau+c_6.\label{eq:H-HP-II-C0-const}
\end{equation}
The constants \(w_1, w_2, w_3\) are given by,
\begin{subequations}
\begin{align}
    w_1&=-\frac{1+\beta}{\gamma}\hat{\rho},\\
    w_2&=-\frac{2\sigma \hat{\rho}}{\gamma \sigma+(1+\beta) \hat{\rho}},\\
    w_3&=\frac{
    (1+\beta)(\gamma\sigma+\hat{\rho}(1+\beta))-2\gamma\sigma
    }
    {2\gamma}\hat{\rho},
\end{align}
\end{subequations}
and the constant of integration \(c_6\) is found by matching to the region I asymptote \(1-\sigma \phi\), giving the value \(c_6 = 1-\sigma \phi - w_2\).

Substituting equations~\eqref{eq:H-HP-II-Q0} and \eqref{eq:H-HP-II-C0-const} into equation~\eqref{eq:H-HP-II-leading-d} and rearranging for \(\tilde{I}_0\) then completes the leading order solution in region II. This expression is quite lengthy to write out in full so the solution is expressed below in terms of \(Q_0\) and \(C_0\). The leading order solutions in region II are therefore,
\begin{subequations}
    \begin{align}
        P_0(\tau) & = 1, \\
        Q_0(\tau) & = -\frac{b}{\gamma \sigma}\left[\frac{2b\tanh\left(b\tau\right)+a}{2b+a\tanh\left(b\tau\right)}\right]+\frac{a}{2\gamma \sigma}, \\
        C_0(\tau) & = w_1\ln{\bigg\lvert \cosh(b\tau)+\frac{a}{2b}\sinh(b\tau)\bigg\rvert}-w_{2}\frac{a\tanh(b\tau)}{2b+a\tanh\left(b\tau\right)}+w_{3}\tau + 1 - \sigma\phi, \\
        \tilde{I}_0(\tau) & = \frac{\gamma\sigma (1-Q_0(\tau))Q_0(\tau)}{C_0(\tau)}.
    \end{align}
\end{subequations}


The switchover time can again be approximated by solving \(C_0(\tau_{sw})\approx 0\). While one could in principle solve this equation numerically for \(\tau_{sw}\), we will instead carry out a series of approximations to yield a closed-form expression. First, we assume that \(b\tau_{sw}\) is sufficiently large to make the approximations \(\sinh(b\tau_{sw})\approx \cosh(b \tau_{sw})\approx \exp(b\tau_{sw})\) and \(\tanh(b\tau_{sw})\approx 1\), which leads to the equation,
\begin{equation}
    0=C_0(\tau_{sw})\approx \left(w_1b + w_3\right)\tau_{sw}+ w_1\ln{\left(\frac{1}{2}+\frac{a}{4b}\right)}-w_2\frac{a}{\left(2b+a\right)}+1-\sigma\phi.\label{eq:H-HP-C0-exp-approx}
\end{equation}
Therefore,
\begin{equation}
    \tau_{sw}\approx \left(-w_1\ln{\left(\frac{1}{2}+\frac{a}{4b}\right)}+ w_2\frac{a}{2b+a}-1+\sigma\phi\right)\left(w_1b+w_3\right)^{-1}. \label{eq:H-HP-tausw}
\end{equation}

Equation~\eqref{eq:H-HP-tausw} follows from a controlled approximation applied to the leading order solution, and could in principle be used for fitting experimental data. However we will instead attempt to simplify this expression further to reduce the number of free parameters, and also to make contact with the result derived by Kerr et al.~\cite{kerr2019mathematical} under the assumptions of quadratic kinetics for the slow reaction. 

The numerator of equation~\eqref{eq:H-HP-tausw} can be written in terms of the dimensionless parameter groupings and \(a\) as,
\begin{equation}
    \left[\frac{\hat{\rho}(1+\beta)}{\gamma}\ln\left(\frac{1}{2}+\frac{a}{2\sqrt{a^2-4\gamma\sigma\hat{\rho}}}\right) - \frac{2\sigma\hat{\rho}}{\sqrt{a^2-4\gamma\sigma\hat{\rho}}+a}\right] - 1 + \sigma\phi.
\end{equation}
As mentioned above, \(a^2>4\gamma\sigma\hat{\rho}\), therefore the logarithmic and square root terms may be expanded in powers of \(4\gamma\sigma\hat{\rho}/a^2\), which gives the following leading order approximation for the numerator in equation~\eqref{eq:H-HP-tausw},
\begin{align}
    & \left[\frac{\hat{\rho}(1+\beta)}{\gamma}\left(\frac{\gamma\sigma\hat{\rho}}{a^2} + O\left(\frac{\gamma\sigma\hat{\rho}}{a^2}\right)^2\right) - \frac{\sigma\hat{\rho}}{a} \left(1+O\left(\frac{\gamma\sigma\hat{\rho}}{a^2}\right)^2\right)\right] - 1 + \sigma\phi, \nonumber \\
    & = \left[-\frac{\gamma}{(\hat{\rho}(1+\beta)/\sigma+\gamma)^2}\hat{\rho} + O\left(\frac{\gamma\sigma\hat{\rho}}{a^2}\right)^2 \right] - 1 + \sigma\phi. \label{eq:H-HP-tausw-exp}
\end{align}
We will now consider conditions under which the term in square brackets in equation~\eqref{eq:H-HP-tausw-exp} can be neglected. Considering all positive values for \(\gamma\), the leading order term is bounded by its value when \(\gamma=\hat{\rho}(1+\beta)/\sigma\), in which case the term takes the value \(\sigma/(4(1+\beta))\). Therefore, provided that \(\sigma\) is small in comparison with \(4(1+\beta)\), the numerator of equation~\eqref{eq:H-HP-tausw} can then be approximated by \(-1+\sigma\phi\). In dimensional parameters, this assumption corresponds to \(n_0/c_0 \ll 4(1+k_4/k_2)\), and so will be valid provided that the initial vitamin C concentration is chosen sufficiently high. 

Under this assumption, the dimensionless expression then simplifies to,
\begin{equation}
    \tau_{sw}\approx \frac{2\gamma(1-\sigma\phi)}{\hat{\rho}\left[(1+\beta)\left(\sqrt{a^2-4\gamma\sigma\hat{\rho}}-a\right)+2\gamma\sigma\right]} \label{eq:H-HP-tsw-simplfied1}
\end{equation}
or in dimensional variables,
\begin{align}
    & t_{sw} = \frac{2k_3(c_0-n_0 \phi)}{n_0 p_0} \nonumber \\
    & 
    \cdot \left([k_2+k_4]\left(\left[\left([k_2+k_4] p_0/n_0+k_3\right)^2-4k_2k_3 p_0/n_0\right]^{1/2}
    -[k_2+k_4] p_0/n_0-k_3\right)+2k_2 k_3\right)^{-1}\mkern-24mu .\label{eq:H-HP-tsw-dim-simplified2}
\end{align}
Equation~\eqref{eq:H-HP-tsw-dim-simplified2} can be considered a somewhat simplified expression for the switchover time, the simplifications being transparently (if not unconditionally) justified. Because there are four free parameters (\(\phi,k_2,k_3,k_4\)) we will find it convenient to make one further assumption, which is that the denominator is dominated by the \(2k_2 k_3\) term. This assumption is uncontrolled and cannot be justified a priori, but has the considerable advantage of reducing the expression to the simple form,
\begin{equation}
    t_{sw}^{\text{H-HP}} := \frac{c_0-\phi n_0}{k_2 n_0 p_0}.\label{eq:H-HP-tsw-simplified}
\end{equation}
Equation~\eqref{eq:H-HP-tsw-simplified} has a similar structure to that of Kerr et al.\ \cite{kerr2019mathematical} and like equation~\eqref{eq:M-HP-tsw} only possesses two fitting parameters, \(k_2\) and \(\phi\). The reasonableness of this step will be demonstrated through simultaneous fitting to experimental data in section~\ref{sec:exp}.

Because equation~\eqref{eq:H-HP-tsw-dim-simplified2} includes an inverse dependence on \(p_0\), the question may be posed (see open peer review accompanying this article) whether this dependence predicts that switchover time can be made arbitrarily small in practice simply by increasing \(p_0\)? However, sufficiently large values of \(p_0\) would correspond to \(p_0/c_0\) no longer being merely \(O(\epsilon^{-1})\), which would be outside of the range of validity of the present analysis by changing the relative orders of magnitudes of the reaction steps. The implications of taking \(p_0/c_0\) at the next order of magnitude \(O(\epsilon^{-2})\) are discussed briefly in Appendix~\ref{sec:AppVHP}; in short, the analysis predicts that the switchover time approaches a non-zero value which has no quantitative dependence on \(p_0\).

\subsection{Region III: corner}
To connect the induction period to the long term state it is again necessary to consider a corner region centred at \(t_{sw}\). The most structured balance of the high peroxide system occurs with timescale \(\epsilon^{1/2}\bar{\tau}=\epsilon t - \tau_{sw}\) and dependent variables scaled as \(C(t)=\epsilon^{1/2}\bar{C}(\bar{\tau})\) and \(I(t)=\epsilon^{1/2} \bar{I}(\bar{\tau})\), with \(Q\) and \(P\) remaining order 1. The scaled system is then,
\begin{subequations}
\begin{align}
    \epsilon^{1/2}\frac{dP}{d\bar{\tau}} & = -\epsilon^2 \sigma(1-Q-2\epsilon^{1/2}\bar{I})P - \epsilon^2 \beta\sigma QP,\\
    \epsilon^{1/2}\frac{dQ}{d\bar{\tau}} & = \epsilon\hat{\rho}(1-Q-2\epsilon^{1/2}\bar{I})P-\epsilon\gamma\sigma(1-Q-2\epsilon^{1/2}\bar{I})Q-\epsilon\beta\hat{\rho}QP,\\
    \frac{d\bar{C}}{d\bar{\tau}} & = -\sigma\bar{I}\bar{C},\\
    \frac{d\bar{I}}{d\bar{\tau}} & = \gamma\sigma(1-Q-2\epsilon^{1/2}\bar{I})Q- \bar{I}\bar{C}.
\end{align}
\end{subequations}
Expanding in powers of \(\epsilon^{1/2}\), the leading order system takes the form,
\begin{subequations}
\begin{align}
    \frac{dP_0}{d\bar{\tau}}       & = 0, \\
    \frac{dQ_0}{d\bar{\tau}}       & = 0, \\
    \frac{d\bar{C}_0}{d\bar{\tau}} & = -\sigma\bar{I}_0\bar{C}_0, \\
    \frac{d\bar{I}_0}{d\bar{\tau}}  & = \gamma\sigma(1-Q_0)Q_0 - \bar{I}_0\bar{C}_0.
\end{align}
\end{subequations}

The variables \(P_0\) and \(Q_0\) are therefore constant in region III, their values being determined by matching to region II. As in region II, \(P_0(\bar{\tau})=1\). From inspection of equation~\eqref{eq:H-HP-II-Q0} and the properties of the \(\tanh\) function, it is clear that the region II solution for \(Q_0\) rapidly approaches the steady state value,
\begin{equation}
   Q_0= Q_0^{ss}=-\frac{b}{\gamma\sigma}\left[\frac{2b+a}{2b+a}\right]+\frac{a}{2\gamma\sigma}.\label{eq:H-HP-III-Q0ss}
\end{equation}
For brevity we will define \(\mu:=\sigma\sqrt{\gamma(1-Q_0^{ss})Q_0^{ss}};\) the remaining equations are,
\begin{subequations}
\begin{align}
    \frac{d\bar{C}_0}{d\bar{\tau}} & = -\sigma\bar{I}_0\bar{C}_0, \\
    \frac{d\bar{I}_0}{d\bar{\tau}} & = \frac{\mu^2}{\sigma}-\bar{I}_0\bar{C}_0.
\end{align}
\end{subequations}
As in section~\ref{sec:M-HP}\ref{sec:M-HP-III}, the variable \(\bar{I}_0\) can be eliminated to yield a Ricatti equation,
\begin{equation}
    \frac{d\bar{C}_0}{d\bar{\tau}} = -(\mu^2\bar{\tau}+c_7+\bar{C}_0)\bar{C}_0,
\end{equation}
where \(c_7\) is a constant of integration, which is again indeterminate at leading order, so it is taken to be zero. This equation takes the same form as that given in section~\ref{sec:M-HP}\ref{sec:M-HP-III} and therefore the solution has a similar form, and \(I_0(\bar{\tau})\) can then be deduced from the relation \(I_0=(C_0+\mu^2\bar{\tau})/\sigma\).

The leading order solutions in region III are therefore,
\begin{subequations}
    \begin{align}
        P_0(\bar{\tau}) & = 1, \\
        Q_0(\bar{\tau}) & = \frac{\left(\hat{\rho}+\gamma \sigma+\beta \hat{\rho}\right)-\sqrt{\left(\hat{\rho}+\gamma \sigma+\beta \hat{\rho}\right)^2-4\gamma \sigma\hat{\rho}}}{2\gamma\sigma}, \\
        \bar{C}_0(\bar{\tau}) & = \frac{\mu\sqrt{2}\exp\left[-\mu^2\bar{\tau}^2/2\right]}{\sqrt{\pi}\left[\mathrm{erf}\left(\mu\bar{\tau}/\sqrt{2}\right)+1\right]}, \\
        \bar{I}_0(\bar{\tau}) & = \frac{\mu\sqrt{2}\exp\left[-\mu^2\bar{\tau}^2/2\right]}{\sigma\sqrt{\pi}\left[\mathrm{erf}\left(\mu\bar{\tau}/\sqrt{2}\right)+1\right]} + \frac{\mu^2}{\sigma} \bar{\tau},
    \end{align}
\end{subequations}
with \(\mu\) as defined in equation~\eqref{eq:H-HP-III-Q0ss} \emph{et seq.}

\subsection{Region IV: equilibration}
As for the M-HP case, following the switchover, \(C\) is again at most order \(\epsilon\) and \(I\), \(P\) are order 1. By contrast with the M-HP case, \(Q\) is also order 1, and the time variable is order \(\epsilon^{-1}\). Setting \(\tilde{\tau}=\epsilon t\) and \(\epsilon\tilde{C}(\tilde{\tau})=C(t)\), the system describing the long-term dynamics after the switchover is then,
\begin{subequations}
\begin{align}
    \frac{dP}{d\tilde{\tau}} & = -\epsilon\sigma(1-Q-2I)P-\epsilon\beta\sigma Q P, \\
    \frac{dQ}{d\tilde{\tau}} & = \hat{\rho}(1-Q-2I)P-\gamma\sigma (1-Q-2I)Q - \beta\hat{\rho} Q P, \\
    \epsilon\frac{dC}{d\tilde{\tau}} & = - \sigma I \tilde{C}, \\
    \frac{dI}{d\tilde{\tau}} & = \gamma \sigma(1-Q-2I) Q - I\tilde{C}.
\end{align}
\end{subequations}

The leading order system is,
\begin{subequations}
\begin{align}
    \frac{dP_0}{d\tilde{\tau}} & = 0, \\
    \frac{dQ_0}{d\tilde{\tau}} & = \hat{\rho}(1-Q_0-2I_0)P_0 - \gamma\sigma(1-Q_0-2I_0)Q_0 - \beta\hat{\rho}Q_0 P_0, \\
    0 & = -\sigma I_0 \tilde{C}_0, \\
    \frac{dI_0}{d\tilde{\tau}} & = \gamma\sigma(1-Q_0-2I_0)Q_0 - I_0\tilde{C}_0.
\end{align}
\end{subequations}
Again it is clear that \(P_0(\tilde{\tau})=1\) and that \(\tilde{C}_0(\tilde{\tau})=0\); indeed \(\tilde{C}_n(\tilde{\tau})=0\) for all \(n\geqslant 0\). The dynamics of \(Q_0\) and \(I_0\) are therefore described by the pair of equations,
\begin{subequations}
\begin{align}
    \frac{dQ_0}{d\tilde{\tau}} & = \hat{\rho}(1-Q_0-2I_0) - \gamma\sigma(1-Q_0-2I_0)Q_0 - \beta\hat{\rho}Q_0, \label{eq:H-HP-IV-2d-Q0}\\
    \frac{dI_0}{d\tilde{\tau}} & = \gamma\sigma(1-Q_0-2I_0)Q_0. \label{eq:H-HP-IV-2d-I0}
\end{align}
\end{subequations}

We have not been able to find a closed form solution of the leading order equations~\eqref{eq:H-HP-IV-2d-Q0}--\eqref{eq:H-HP-IV-2d-I0} for order 1 values of the dimensionless parameters. The system has a single equilibrium point \((Q_0,I_0)=(0,1/2)\), which physically corresponds to all iodide and hypoiodous acid having been converted to iodine. The eigenvalues of this equilibrium point are \(0\) and \(-\hat{\rho}(1+\beta)\). The zero eigenvalue corresponds to the slow manifold which is tangent to the line \(\{(2s,1/2-(1+\beta)s \, : \, s \geqslant 0\}\). The negative eigenvalue corresponds to a stable manifold tangent to \(\{(s,0) \, : \, s\in \mathbb{R}\}\) which is outside of the physical domain of interest because it would entail either \(Q\) or \(D\) being negative.

The `terminal dynamics' close to the equilibrium point can be inferred approximately by considering a quasi-steady balance of the terms linear in \(1-Q_0-2I_0\) and \(Q_0\) in equation~\eqref{eq:H-HP-IV-2d-Q0}, yielding,
\begin{equation}
    0\approx \hat{\rho}(1-Q_0-2I_0) - \beta\hat{\rho}Q_0, \label{eq:H-HP-IV-QS}
\end{equation}
which can be rearranged to provide the following approximate relationship between \(I_0\) and \(Q_0\),
\begin{equation}
    I_0 \approx \frac{1-(1+\beta)Q_0}{2},\label{eq:H-HP-IV-QSr}
\end{equation}
from which we deduce that,
\begin{equation}
    \frac{dI_0}{d\tilde{\tau}} \approx -\frac{1+\beta}{2}\frac{dQ_0}{d\tilde{\tau}}. \label{eq:H-HP-IV-QSd}
\end{equation}
Substituting equations~\eqref{eq:H-HP-IV-QS} and \eqref{eq:H-HP-IV-QSd} into equation~\eqref{eq:H-HP-IV-2d-I0} then provides an approximation to the long term dynamics of \(Q_0\):
\begin{equation}
    \frac{dQ_0}{d\tilde{\tau}} \approx -\frac{2\beta\gamma\sigma }{1+\beta} Q_0^2.
\end{equation}
This equation has solution,
\begin{equation}
    Q_0(\tilde{\tau}) = \frac{1+\beta}{2\beta\gamma\sigma \tilde{\tau}+c_8},
\end{equation}
where \(c_8\) is a constant of integration. The corresponding approximation to \(I_0\) then follows from equation~\eqref{eq:H-HP-IV-QSr},
\begin{equation}
    I_0(\tilde{\tau}) \approx \frac{1}{2}\left(1-\frac{(1+\beta)^2}{2\beta\gamma\sigma \tilde{\tau}+c_8}\right).
\end{equation}
The constant \(c_8\) can be determined by matching the solution for \(I_0\) to the region II solution, giving the value \(c_8 = (1+\beta)^2-2\beta\gamma\sigma \tau_{sw}\). The approximate solutions in region IV, expected to be increasingly accurate close to the equilibrium point, are therefore,
\begin{subequations}
\begin{align}
    P_0(\tilde{\tau}) & = 1, \\
    Q_0(\tilde{\tau}) & = \frac{1+\beta}{2\beta\gamma\sigma (\tilde{\tau}-\tau_{sw})+(1+\beta)^2}, \\
    C_0(\tilde{\tau}) & = 0, \\
    I_0(\tilde{\tau}) & = \frac{1}{2}\left(\frac{2\beta\gamma\sigma(\tilde{\tau}-\tau_{sw})}{2\beta\gamma\sigma (\tilde{\tau}-\tau_{sw})+(1+\beta)^2}\right).
\end{align}
\end{subequations}

Figure~\ref{fig:NumericalSolutionLogRegions}(b) compares the high hydrogen peroxide asymptotic solutions in each region with a numerical solution computed with an arbitrary set of parameters (details in caption). The solutions are quite close in their region of validity, particularly for the trajectories of \(P\), \(I\) and \(C\), although as time progresses, a significant drift develops between the numerically computed value of \(P\) and the constant approximation produced by the asymptotic analysis; the discrepancy appears to be on the order of \(\epsilon^{1/2}\). This discrepancy could be reduced by seeking the next order term in the asymptotic expansion, although this quantity is less critical to determine than the clock chemicals \(I\) and \(C\) which directly determine the switchover time.

\begin{figure}
    \centering
    \includegraphics{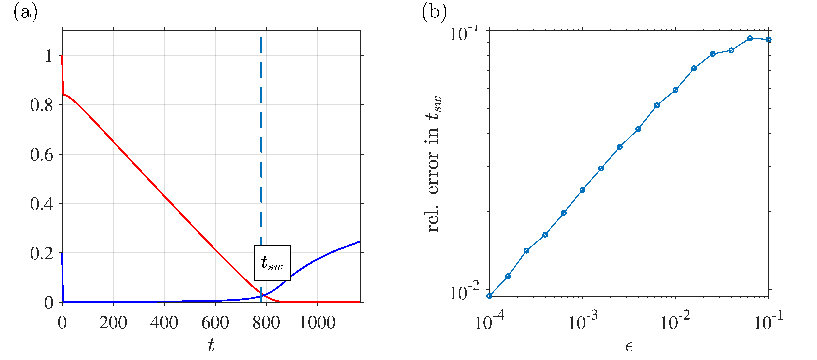}
    \caption{Comparison of asymptotic-derived formula and numerical solution in determination of the dimensionless switchover time for regime H-HP. Parameter values (chosen arbitrarily): \(\epsilon=10^{-2}\), \(\beta=0.6\), \(\gamma=0.7\), \(\sigma=0.8\), \(\phi=0.2\) and \(\hat{\rho}=0.9\). (a) Numerical solutions for \(C(t)\) and \(I(t)\) shown alongside the value of the switchover time (dashed line) calculated by the approximate formula \eqref{eq:H-HP-tausw}. (b) Relative error between the asymptotic expression for the switchover time \eqref{eq:H-HP-tausw} and the point at which the numerical solution falls below a threshold value of \(\epsilon\), as a function of the small parameter \(\epsilon\). }\label{fig:H-HP-Convergence}
\end{figure}

Figure~\ref{fig:H-HP-Convergence}(a) compares the asymptotic-derived switchover time formula \eqref{eq:H-HP-tausw} to that of the numerical solution with the relative error between these for varying \(\epsilon\) given in Figure~\ref{fig:H-HP-Convergence}(b). By contrast with the moderate hydrogen peroxide case, the convergence of the error as \(\epsilon\rightarrow 0\) is sublinear, the gradient on a log-log plot approaching approximately \(0.4\). These results provide evidence that the asymptotic expansion has sublinear accuracy in this region and so a correction with a fractional power in \(\epsilon\) may be needed. We will not pursue this analysis in the present paper as it appears likely to involve significant additional complications for limited additional insight.

Overall, the comparison for each of the M-HP and H-HP problems in figure~\ref{fig:NumericalSolutionLogRegions} provide strong evidence for the validity of the asymptotic approach. In particular, the approximations for the clock chemical iodine are both very satisfactory, and as the results for the inhibitor vitamin C, and the intermediate species hypoiodous acid. For large time (region IV for M-HP and late in region II for H-HP), a discrepancy emerges between the numerical solution and leading order asymptotic solution for hydrogen peroxide, likely due to the need for more terms or possibly a slower-decaying asymptotic series approximation. Because the error in hydrogen peroxide does not appear to affect the approximation to the clock chemical concentration or switchover time, we will not investigate this issue further in the present manuscript.

\section{Experimental testing}\label{sec:exp}
The M-HP and H-HP switchover time formulae \eqref{eq:M-HP-tsw} and \eqref{eq:H-HP-tsw-simplified} were tested through tabletop experiments in a similar manner to the work of Kerr et al.~\cite{kerr2019mathematical}, using a combination of vitamin C powder, Lugol's iodine, hydrogen peroxide solution, and powdered laundry starch. A key refinement to the experimental technique was to use a webcam sensor running under Matlab R2020b to detect the point at which the colour change occurred. All experiments were carried out diluted in 60~ml water at \(40^\circ\), with 5~ml of vitamin C stock solution (made up as 1000~mg in 30--90~ml water), 3--10~ml of 3\% Lugol's iodine, 5~g laundry starch and 1--20~ml 3\% hydrogen peroxide. 

\subsection{Imaging}
The switchover time \(t_{sw}\) was measured by imaging the mixture from above under natural lighting conditions using a RGB USB camera running under Matlab with Image Processing Toolbox (Mathworks, Natick) 2020b, recording images at 15 frames per second. The region of interest was set at \(80\times 100\) pixels in the centre of the beaker. The red channel was sufficient for a clear measurement. Briefly, the signal was processed by summing the intensity across all pixels, taking the difference of successive frames, then taking forward and backward moving averages with a 10 frame window to minimise the effect of fluctuations, then identifying the point at which the difference between forward and backward moving averages took the largest negative value. This `corner' value corresponds approximately to the point of largest second derivative of the signal, and is therefore relatively insensitive to lighting conditions. Image processing code can be found in the accompanying supplement \cite{alsaleh2024}. 

\subsection{Fitting series}
Four experimental series, referred to as \(N_M\) and \(C_M\) (moderate hydrogen peroxide regime), and \(N_H\) and \(C_H\) (high hydrogen peroxide regime) were carried out for fitting of the constants \(\phi\) and \(k_3\). These experiments involved holding varying respectively initial iodine concentration or vitamin C concentration while holding the initial masses of other substances fixed. Due to slightly varying volume of the overall solution, in series \(N_M\) and \(N_H\), the concentrations of vitamin C and hydrogen peroxide also varied slightly, which was taken into account when fitting. Full details of the concentrations is given in the accompanying supplement \cite{alsaleh2024}.

The moderate hydrogen peroxide regime involved using \(1\)~ml hydrogen peroxide, providing an initial concentration of \(p_0\approx 6.7\times 10^{-3}\)~mol/l. In series \(N_M\) the initial mass of Lugol's was varied from \(4\) to \(8.5\)~ml, producing an initial concentration of \(n_0=6.1\times 10^{-3}\) to \(1.2\times 10^{-2}\)~mol/l, with \(c_0\approx 2.3\times 10^{-3}\)~mol/l. In series \(C_M\) the initial mass of vitamin C was varied from \(4.2\times 10^{-2}\) to \(7.1\times 10^{-2}\)~g, producing a concentration of \(c_0=1.8\times 10^{-3}\) to \(7.2\times 10^{-3}\)~mol/l, with \(n_0\approx 7.6\times 10^{-3}\)~mol/l. In all of these experiments, the initial ratio of hydrogen peroxide to iodine was in the range \(0.52\leqslant \rho \leqslant 1.1\). Measured values of the switchover time varied between approximately 60 and 700~s.

The high hydrogen peroxide regime used \(20\)~ml hydrogen peroxide, providing an initial concentration of \(p_0\approx 0.12\)~mol/l. In series \(N_H\) the initial mass of Lugol's was varied from \(3\) to \(10\)~ml, producing an initial concentration of \(n_0=4.0\times 10^{-3}\) to \(1.2\times 10^{-2}\)~mol/l, with \(c_0\approx 6.3\times 10^{-3}\)~mol/l. In series \(C_H\) the initial mass of vitamin C was varied from \(5.6\times 10^{-2}\) to \(0.17\)~g, producing a concentration of \(2.1\times 10^{-3}\) to \(6.3\times 10^{-3}\)~mol/l, with \(n_0\approx 6.6\times 10^{-3}\)~mol/l. In these experiments the ratio \(8.8<p_0/n_0<30\). Measured values of the switchover time varied from 23 to around 190~s.

The values of \(\phi\) and \(k_3\) were estimated through relative error least squares fitting in order to account for the wide range of measured switchover time, carried out in Matlab with the unconstrained local optimisation function \texttt{fminsearch}, which implements the Nelder-Mead simplex search method \cite{nelder1965simplex}. Confidence limits were estimated through bootstrapping with \(N=10000\) samples, implemented in Matlab with the function \texttt{bootci}, which implements the bias corrected and accelerated percentile method \cite{efron1987better}.

\subsection{Testing series}
The switchover time model and fitted parameters were then tested by comparing against two independent data series, referred to as \(H_1\) and \(H_2\). Series \(H_1\) varied the hydrogen peroxide concentration between \(3.3\times 10^{-2}\)~mol/l and \(0.17\)~mol/l, while holding \(c_0\approx 6\times 10^{-3}\)~mol/l and \(n_0\approx 7\times 10^{-3}\)~mol/l, so that \(4.4\leqslant p_0/n_0 \leqslant 27\). Series \(H_2\) varied the hydrogen peroxide concentration between \(6.7\times 10^{-2}\) and \(3.6\times 10^{-2}\)~mol/l, while holding \(c_0\approx 2.4 \times 10^{-3}\)~mol/l and \(n_0\approx 7.5\times 10^{-3}\)~mol/l, so that \(0.88\leqslant p_0/n_0 \leqslant 4.9\). The switchover time formula \eqref{eq:M-HP-tsw} was used for cases where \(p_0/n_0 \leqslant 1.5\) and \eqref{eq:H-HP-tsw-simplified} when \(p_0/n_0>1.5\), the threshold value of \(1.5\) being chosen arbitrarily.

\subsection{Results}
\begin{figure}
    \centering
    \includegraphics{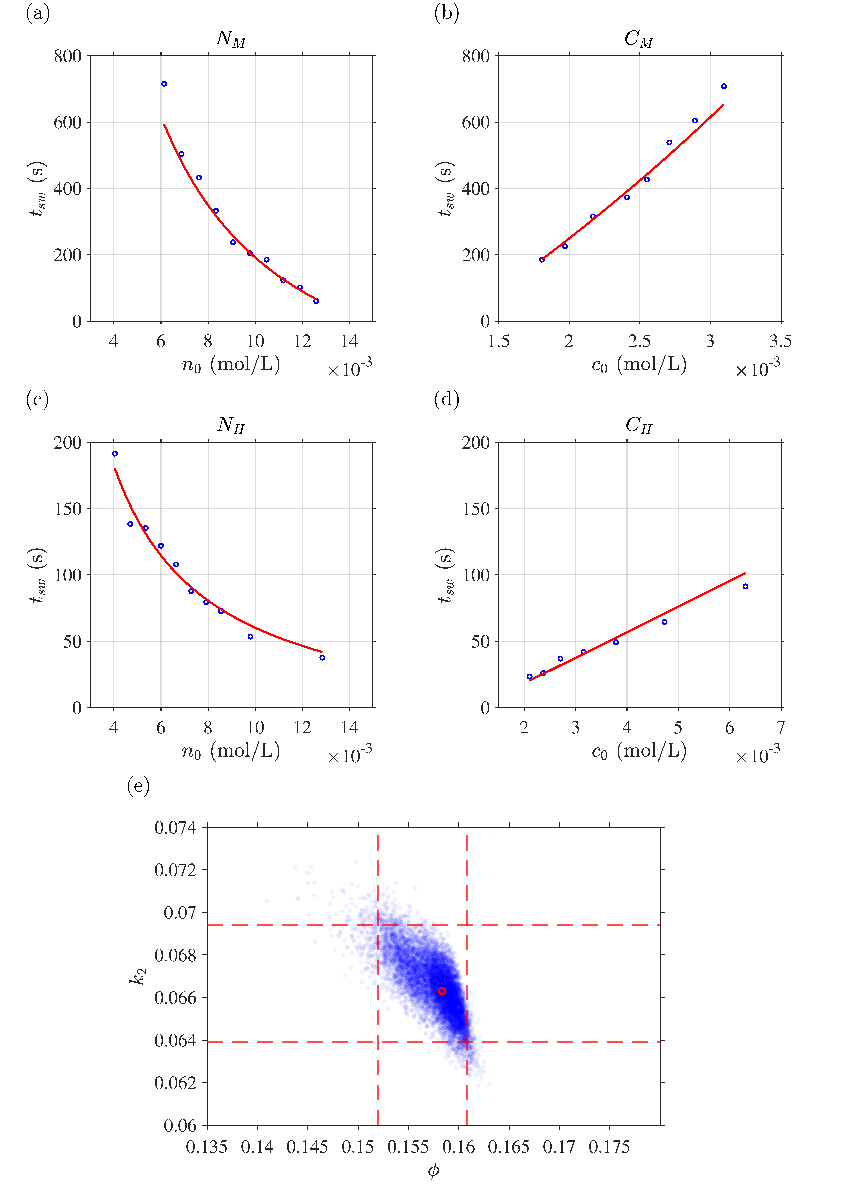}
    \caption{Outcome of simultaneous fitting to data series \(N_M\), \(C_M\), \(N_H\) and \(C_H\). (a--d) experimental (blue dot) and fitted switchover time (red line) with sum square relative error best fit of \(\hat{k}_2=0.066302\) and \(\hat{\phi}=0.15835\). (e) Bootstrapping results with \(10000\) repeats (blue dots), best fit (red circle) and 95\% confidence intervals for each parameter (dashed red lines).} \label{fig:fitting}
\end{figure}

Figure~\ref{fig:fitting}(a--d) shows the outcome of the fitting series experiments \(N_M\), \(C_M\), \(N_H\) and \(C_H\). Qualitatively, the results show that in both low and high hydrogen peroxide regimes, switchover time declines as initial iodine concentration is increased, and increases nearly linearly with initial vitamin C concentration -- fundamentally, as substrate is increased, inhibitor is used up more quickly, whereas as inhibitor is increased, it lasts longer. The quality of the fit appears excellent, with the maximum relative error rarely exceeding 10\%. The estimated parameter values of \(\hat{\phi}=0.158\) (95\% confidence interval \([0.152,0.161]\)) and \(\hat{k}_2=0.0663\) (95\% CI \([0.639,0.693]\) are shown in figure~\ref{fig:fitting}(e) along with bootstrap sample distributions. 

\begin{figure}
    \centering
    \includegraphics{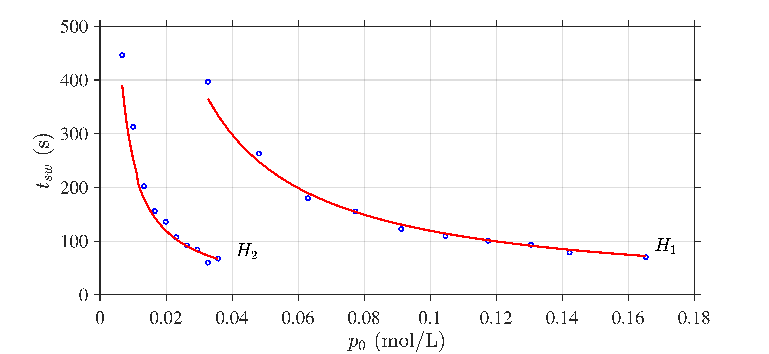}
    \caption{Test of the parameter fits (figure~\ref{fig:fitting}) against independent experimental series \(H_1\) and \(H_2\) in which initial hydrogen peroxide mass is varied.}\label{fig:testing}
\end{figure}

The testing series \(H_1\) and \(H_2\) using the estimated parameters provide at least as close a match between experiment and the fitted model (figure~\ref{fig:testing}) as for the fitted data, providing confidence in the out-of-sample applicability of the model and accuracy of the fitted parameters. Qualitatively, the data show switchover time decreasing as hydrogen peroxide concentration is increased, corresponding to the slow reaction being promoted through additional availability of a key reactant in its rate-limiting step.

\section{Discussion}
This paper unifies and significantly extends previous work on asymptotic analysis and experimental data fitting for a tabletop chemical kinetics experiment. Through breaking down the reaction converting iodine to iodine into two steps involving the intermediate hypoiodous acid, the resulting model was able to explain the kinetics observed in regimes of moderate and high hydrogen peroxide concentration, encompassing both the models contained in Kerr et al.\ \cite{kerr2019mathematical} and Parra Cordova and Gonz\'{a}lez Pe\~{n}a \cite{parra2020enhancing}, and being transparently consistent with previous descriptions of the kinetics of the iodide-hydrogen peroxide reaction \cite{liebhafsky1933kinetics,copper1998kinetics,sattsangi2011microscale}. The key kinetic insight provided by the model is that the production of hypoiodous acid is only rate-limiting in the situation where hydrogen peroxide is comparable to other reactants.

The analysis utilised a single small parameter \(0<\epsilon\ll 1\) which was used to represent the disparity of rates between the slow rate of production of hypoiodous acid from iodide, the production of iodine from iodide and hypoiodous acid, and the conversion of iodine back to iodide through oxidation of vitamin C (fastest). A slow reverse reaction reducing hypoiodous acid back to iodide was also included. The small parameter was thereby able to distinguish the moderate from high hydrogen peroxide regime, the hydrogen peroxide-to-vitamin C ratio being order 1 in the moderate case, and order \(\epsilon^{-1}\) in the high case.

Asymptotic analysis was carried out in each of the two regimes. The moderate hydrogen peroxide case involved five distinct time regions: an \(O(1)\) timescale (region I) during which the initial free iodine is mostly converted to iodide through reduction of a fraction of the initial vitamin C. This region is followed by an \(O(\epsilon^{-1})\) process (region Ia) during which iodide, iodine and hypoiodous acid reach quasi equilibrium. When comparing against numerical solutions, the leading order terms captured the kinetics of these species well, although there was some divergence from the hydrogen peroxide and vitamin C concentrations in this region. The induction period (region II) followed, which occurred on \(O(\epsilon^{-2})\) time. This process captured the consumption of hydrogen peroxide and vitamin C quite accurately, and provided a closed form solution \eqref{eq:M-HP-tsw} for the switchover time, through solving for the leading order solution for the vitamin C concentration reaching zero. This formula 
can be seen as a refinement to that given by \cite{parra2020enhancing} by taking into account the effect of initially non-zero molecular iodine in oxidising vitamin C and hence shortening the induction period. By similar methods to Kerr et al.\ it was also possible to construct an approximate solution in the intermediate region (III) for order \(\epsilon^{-1}\) time around the switchover, and the long term evolution to equilibrium (region IV) on a timescale of order \(\epsilon^{-2}\). The leading order solutions on regions II--IV showed reasonably good agreement with numerical solutions, in particular the relative error between the switchover time formula and the point at which the numerical solution fell below \(\epsilon\) was found to converge approximately linearly as \(\epsilon\rightarrow 0\). An area for possible future development would be to predict better the consumption of vitamin C and hydrogen peroxide in region Ia, and to provide a higher order approximation in region II, however as noted above, it seems unlikely that significant additional insight would be obtained.

The second asymptotic analysis of the paper concerned the high hydrogen peroxide case already considered by Kerr et al. Heuristically one can argue that when hydrogen peroxide levels are sufficiently high, the rate-limiting step is no longer such, and so the overall reaction \(2I^{-}\rightarrow I_2\) can be modelled by the quadratic kinetics that would follow from the law of mass action, and hydrogen peroxide levels need not be modelled explicitly as they are not significantly reduced. The analysis of section~\ref{sec:H-HP} put this informal argument on a more quantitative basis through elucidating how a hydrogen peroxide concentration of order \(\epsilon^{-1}\) provides the balance described above. Similarly to Kerr et al.\ we found four asymptotic regions corresponding to (I) initial reduction of vitamin C due to non-zero initial iodine, which is simultaneous with iodide, iodine and hypoiodous acid reaching quasi-equilibrium, (II) induction, (III) crossover and (IV) evolution to final equilibrium. The leading order solution for region IV could not be solved exactly, however a linearisation about the equilibrium combined with matching to region II enabled a relatively accurate approximate solution to be found; however the leading order analysis was not even in principle able to capture the significant consumption of hydrogen peroxide seen in region IV. 

One of the goals of this paper has been to show how previously-reported formulae for the switchover time can be recovered from our unified model, when considered in the appropriate concentration regime. As previously, the region II solution enabled the deduction of a `full' switchover time formula \eqref{eq:H-HP-tausw} for the full H-HP regime. This full formula was shown to be equivalent to a simplified version that could be justified rigorously under certain assumptions regarding parameters \eqref{eq:H-HP-tsw-dim-simplified2} and a highly simplified heuristic version for comparison with Kerr et al.\ \eqref{eq:H-HP-tsw-simplified}. Due to involving only two free parameters, it was the latter formula that we subsequently used for fitting experimental data. The full formula was found to converge to the numerical solution as \(\epsilon\rightarrow 0\) sublinearly, suggesting that the next order correction may be sublinear. 

The switchover time formulae were then tested through fitting the parameters \(\phi\) (dimensionless initial proportion of molecular iodine) and \(k_2\) (production of hypoiodous acid from iodide and hydrogen peroxide). Parameters were fitted through minimising least square relative error in order to provide similar weighting to results across more than an order of magnitude in time values. Simultaneous fitting was carried out by varying the initial mass of iodine or vitamin C, repeated in moderate and high hydrogen peroxide regimes -- four series to estimate two parameters. Uncertainties were quantified through bootstrapping, with the 95\% confidence regions for each parameter being on the order of \(10\)\% of the estimated value. The formula and estimates were then assessed against two separate experimental series in which initial mass of hydrogen peroxide was varied from moderate to high values. The model curves fell very close to the experimental results without any further fitting or tuning, giving high confidence in both the model and estimates.

The paper thereby presents a validated mechanistic understanding of a paradigm chemical kinetics system. We believe the work may form a useful basis to develop quantitive understanding of other substrate-depletive clock reactions and their applications, in addition to providing an accessible case study in the combination of time-dependent asymptotic analysis and mathematical model parameterisation, an approach which is increasingly used to study biological and biochemical systems \cite{youlden2021}. There are several areas in which the study may be developed and improved further. 

In respect of the `mathematical' aspects, in matching the induction region to the corner region, it was noted that more accuracy could have been obtained if a higher order approximation to the induction region solution were available; this issue may also underlie the sublinear convergence of the numerical and leading order solutions for the switchover time observed in the high hydrogen peroxide model. It was also noticeable that the leading order constant solution for hydrogen peroxide loses accuracy at large time values; it would be of interest to explore whether these dynamics could easily be captured through a higher order asymptotic solution.

Focusing on data fitting, while the strategy taken in this paper was to carry out data fitting on a highly simplified form of the high hydrogen peroxide switchover time formula, it should in principle be possible to fit the full form following from equation~\eqref{eq:H-HP-tausw}. The fitted model parameters may also provide post hoc justification of the simplifications carried out to form equation~\eqref{eq:H-HP-tsw-dim-simplified2}. 

From a chemical kinetics modelling perspective, the hypoiodous acid-producing reaction has been found to comprise two parallel pathways, one of which depends on the availability of \(H^+\) ions \cite{copper1998kinetics}. Therefore it would be of interest to model these pathways independently in order to predict the effect of pH on the overall reaction rate. While the present study fixed the temperature of the reaction, another area of development would be to model temperature dependence of reaction rates and fit the model parameters by varying temperature systematically. 

In respect of experimental technique, another area for development would be the use of ultraviolet light spectrophotometry \cite{sulistyarti2015} to measure the time evolution of iodine levels (via blue iodine-starch complex) rather than simply the switchover time, so that model curves could be fitted directly, providing a more complete test of the accuracy of the solutions in each asymptotic region. Explicit mathematical modelling of the formation of starch-iodine complex may in turn be required to interpret spectrophotometry data. As the discussion above illustrates, clock reactions continue to provide challenge and inspiration at the interface of mathematics and physical chemistry.

\vspace{\baselineskip}

{
\noindent {\bf Data accessibility:} {data are available at the Dryad repository \cite{alsaleh2024}}. \\
\noindent {\bf Competing interests:} we declare we have no competing interests. \\
\noindent {\bf Authors' contributions:} AAS formulated the model, conducted asymptotic analysis and experiments and wrote computer code. DJS and SJ designed the research and provided supervision and support. All authors contributed to analysing results and drafting/editing the manuscript. \\
\noindent {\bf Funding statement:} AAS acknowledges support from the Kingdom of Saudi Arabia, Ministry of Education, The Custodian of the Two Holy Mosques Scholarship Program (12th), KSP12004945, Ref. 1038615611.
}

\appendix
\section{Quasi-equilibrium for \(Q(t)\) in the M-HP regime}\label{sec:appendix}
\begin{figure}
    \centering
    \includegraphics{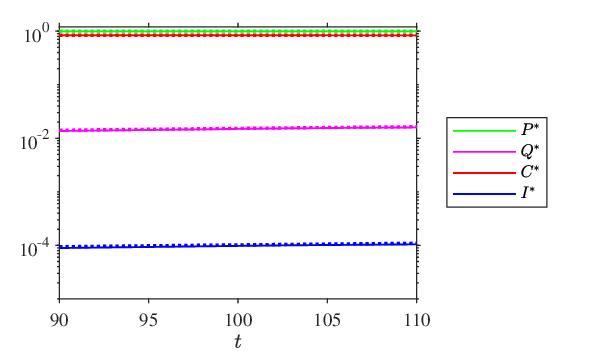}
    \caption{Comparison of numerical solutions (solid line) and asymptotic solutions (dotted line), replotted from figure~\ref{fig:NumericalSolutionLogRegions}(a), region Ia to resolve the dynamics for \(Q(t)\).}\label{fig:Comparison-M-Ia}
\end{figure}

The establishment of quasi-equilibrium for \(Q(t)\) (Figure 1(a), region Ia) is displayed in more detail in figure~\ref{fig:Comparison-M-Ia}, with vertical axis for concentration and linear axis for time. It is evident that \(Q(t)=O(\epsilon)\) and \(I(t)=O(\epsilon^2)\).

\section{Very high hydrogen peroxide regime}\label{sec:AppVHP}
During open peer review, a query was raised as to the behaviour of \(t_{sw}\) as \(p_0\) is made arbitrarily large -- specifically, does it follow from equation~\eqref{eq:H-HP-tsw-dim-simplified2} that \(t_{sw}\) should approach zero? However, the results of section~\ref{sec:H-HP} rest on the assumption that the initial hydrogen peroxide level is on the order of \(\epsilon^{-1}\) relative to vitamin C, i.e.\ \(\rho=\epsilon^{-1}\hat{\rho}\) with \(\hat{\rho}\) being order 1. To explore the next order of magnitude in relative hydrogen peroxide concentration, below we briefly consider how the dynamics up to switchover change when instead \(\rho=\epsilon^{-2}\hat{\rho}\). The dimensionless equation for \(Q\) is then,
\begin{equation}
    \frac{dQ}{dt} = \hat{\rho} (1-Q-2I)P - \epsilon\gamma\sigma(1-Q-2I)Q - \beta \hat{\rho} QP.
\end{equation}
Qualitatively, the first and third terms on the right hand side are now order 1, which means that the ``\(k_2\)'' and ``\(k_4\)'' reactions predominate, and \(Q\) reaches an order 1 quasi-equilibrium in order 1 time. In region I, the leading order behaviour is then,
\begin{equation}
    \frac{dQ_0}{dt} = \hat{\rho} \left(1-(1+\beta)Q_0-2I_0\right) P_0,
\end{equation}
which approaches a steady state \(Q_0(t)\rightarrow (1+\beta)^{-1}\) as \(t\rightarrow \infty\).

Region II then again corresponds to the rescaling \(t=\epsilon^{-1}\tau\) and \(I=\epsilon\tilde{I}\); the dynamics for \(Q(\tau)\) are then,
\begin{equation}
    \frac{dQ}{d\tau} = \epsilon^{-1}\hat{\rho}\left(1-Q-2\epsilon\tilde{I}\right)P - \gamma \sigma\left(1-Q-2\epsilon\tilde{I}\right) Q - \epsilon^{-1} \beta\hat{\rho}Q P,
\end{equation}
which at leading order yields the constant solution \(Q_0(\tau)=(1+\beta)^{-1}\). The remaining leading order equations for \(P_0\), \(\tilde{I}_0\) and \(C_0\) are unchanged from equations~\eqref{eq:H-HP-II-leading}. However with the modified value of \(Q_0\), the solution for \(C_0\) then becomes,
\begin{align}
    C_0(\tau) = 1 - \sigma\phi - \frac{\beta\gamma\sigma^2}{(1+\beta)^2}\tau,
\end{align}
from which it follows that,
\begin{equation}
    \tau_{sw}=\frac{(1-\sigma\phi)(1+\beta)^2)}{\beta\gamma\sigma^2}.
\end{equation}
Returning to dimensional variables, the switchover time formula in this regime is then,
\begin{equation}
    t_{sw} = \frac{c_0-\phi n_0}{(k_3/k_1)(k_4/k_2) (1+(k_4/k_2))^{-2} n_0^2},
\end{equation}
which takes a similar form to the formula of Kerr et al.\ \cite{kerr2019mathematical}. The switchover time depends on the reaction rates only through the ratios \(k_3/k_1\) (2nd forward step of slow reaction relative to fast reaction) and \(k_4/k_2\) (reverse step of slow reaction relative to first forward step). Qualitatively, it is the second (\(k_3\)) step that is now rate limiting for the slow reaction, with the rate set by the \(Q\)-level that emerges from the quasi-steady balance of the hydrogen peroxide-dependent reactions (\(k_2\) and \(k_4\)).

In conclusion, we predict that once the hydrogen peroxide level is sufficiently high, the analysis of section~\ref{sec:H-HP} is no longer valid and \(t_{sw}\) does not approach zero. Instead, the switchover time formula depends on vitamin C and iodine concentrations similarly to the high hydrogen peroxide case, but the quantitative dependence on hydrogen peroxide level is lost.

{
\begin{multicols}{2}
    \raggedright

\end{multicols}
}


\end{document}